\pdfoutput=1
\RequirePackage{silence}
\documentclass[a4paper,svgnames]{amsart}
\usepackage[hmarginratio={1:1},vmarginratio={1:1},lmargin=60.0pt,tmargin=60.0pt]{geometry}

\usepackage{booktabs}
\allowdisplaybreaks

\usepackage[T1]{fontenc}
\usepackage{latexsym,exscale,amsfonts,amssymb,mathtools}
\usepackage{amsmath,amsthm,amsfonts,amssymb,amscd,textcomp,bbm}
\usepackage{mathrsfs,stackrel}
\usepackage{etoolbox}
\usepackage{tikz,tikz-cd}
\usetikzlibrary{matrix,arrows.meta,decorations.markings,shapes.multipart,decorations.pathreplacing,decorations.pathmorphing}
\usepackage{aliascnt}
\usepackage{ytableau}
\usepackage{xparse}
\usepackage{cite}

\usepackage{dynkin-diagrams}

\usepackage{array}
\newcolumntype{C}{>{$}c<{$}}

\usepackage{xcolor,colortbl}

\definecolor{mygray}{gray}{0.6}
\definecolor{mygraydark}{gray}{0.4}
\definecolor{mygraylight}{gray}{0.85}
\definecolor{spinach}{RGB}{46,139,87}
\definecolor{tomato}{RGB}{255,99,71}
\definecolor{orchid}{RGB}{143,40,194}
\definecolor{neon}{RGB}{77,77,255}
\definecolor{pumpkin}{RGB}{224,180,80}
\definecolor{citron}{RGB}{190,180,90}

\definecolor{lava}{RGB}{207,16,32}
\definecolor{cream}{RGB}{255,253,208}
\definecolor{verdigris}{RGB}{67,179,174}
\definecolor{Black}{RGB}{0,0,0}
\definecolor{mydarkblue}{RGB}{10,10,170}
\definecolor{darkspinach}{RGB}{20,70,20}
\definecolor{darktomato}{RGB}{155,40,30}
\definecolor{darkorchid}{RGB}{50,10,100}
\definecolor{darklava}{RGB}{150,8,16}

\usepackage{todonotes}

\def\changed#1{{#1}}
\def\ochanged#1{{#1}}
\def\tchanged#1{{#1}}

\usepackage{enumitem}
\setlist[enumerate]{itemsep=0.15cm,label=\emph{\upshape(\alph*)}}
\setlist[enumerate,2]{itemsep=0.15cm,label=\emph{\upshape(\roman*)}}
\setlist[enumerate,3]{itemsep=0.15cm,label=\emph{\upshape(\Alph*)}}

\let\emph\relax
\DeclareTextFontCommand{\emph}{\em}

\renewcommand{\dots}{\text{...}}

\newcommand{\rddots}{\raisebox{-0.075cm}{\rotatebox{45}{\text{...}}}}
\newcommand{\placeholder}{{}_{-}}

\newcommand{\cf}{\text{cf.}}

\newcommand{\C}{\mathbb{C}}

\newcommand{\N}{\mathbb{Z}_{\geq 0}}

\newcommand{\Z}{\mathbb{Z}}
\newcommand{\K}{\mathbb{K}}

\newcommand{\vcirc}{\circ}
\newcommand{\hcirc}{\otimes}
\newcommand{\munit}{\mathbbm{1}}

\newcommand{\setstuff}[1]{\mathrm{#1}}
\newcommand{\catstuff}[1]{\mathbf{#1}}
\newcommand{\functorstuff}[1]{\mathcal{#1}}

\font\scc=rsfs10
\newcommand{\twocatstuff}[1]{\scc\mbox{#1}\hspace{1.0pt}}

\newcommand{\obstuff}[1]{\mathtt{#1}}
\newcommand{\morstuff}[1]{\mathrm{#1}}

\newcommand{\idmor}{\morstuff{id}}

\newcommand{\End}{\setstuff{End}}
\newcommand{\Hom}{\setstuff{Hom}}

\newcommand{\bx}[1][x]{\vec{#1}}
\newcommand{\bN}[1][N]{\vec{#1}}

\newcommand{\matall}[1][n]{\mathrm{Mat}_{#1}(\C)}
\newcommand{\matgln}[1][n]{\mathrm{GL}_{#1}(\C)}
\newcommand{\matsln}[1][2]{\mathrm{SL}_{#1}(\C)}
\newcommand{\mat}[1][N]{\mathsf{#1}}
\newcommand{\trace}{\mathrm{tr}}

\newcommand{\congruence}{\equiv_{c}}

\newcommand{\tensors}[3]{\mathrm{T}_{#1,#2,#3}(\C)}

\newcommand{\ten}[1][T]{\mat[#1]}

\newcommand{\trep}[1][M]{\functorstuff{#1}}
\newcommand{\mrep}{\functorstuff{F}}

\newcommand{\equirep}{\cong_{2rep}}

\newcommand{\fincat}[1][\K]{\twocatstuff{A}}
\newcommand{\ffincat}[1][\K]{\scalebox{0.85}{$\twocatstuff{A}$}}
\newcommand{\rank}{\mathrm{rank}\,}

\newcommand{\vect}[1][\C]{#1\catstuff{Vect}}
\newcommand{\rep}[1][\C]{#1\catstuff{Rep}}
\newcommand{\corep}[1][\C]{#1\catstuff{coRep}}

\newcommand{\tworep}[1][M]{\twocatstuff{#1}}

\newcommand{\dsum}{\subset_{\oplus}}
\newcommand{\add}[1]{{#1}^{\oplus}}
\newcommand{\addc}[1]{{#1}^{\oplus,\dsum}}

\newcommand{\addg}[1][M]{K_{0}^{\oplus}(#1)}
\newcommand{\addgg}[1][M]{K_{0}^{\oplus}\big(#1\big)}

\newcommand{\equirept}{\cong_{rep}}

\newcommand{\sln}[1][2]{\mathrm{SL}_{#1}}
\newcommand{\son}[1][3]{\mathrm{SO}_{#1}}
\newcommand{\gln}[1][2]{\mathrm{GL}_{#1}}

\newcommand{\qpar}{q}
\newcommand{\qnum}[1][2]{[#1]_{\qpar}}

\newcommand{\xpar}{x}

\newcommand{\web}[1][{\sln[2]}]{\twocatstuff{W}\catstuff{eb}(#1)}
\newcommand{\wweb}[1][{\sln[2]}]{\scalebox{0.8}{$\twocatstuff{W}\catstuff{eb}(#1)$}}
\newcommand{\genob}[1][X]{\obstuff{#1}}

\newcommand{\peb}[1][{\gln[2]}]{\twocatstuff{P}\catstuff{eb}(#1)}
\newcommand{\qeb}[1][{\gln[2]}]{\twocatstuff{U}\catstuff{eb}(#1)}

\newcommand{\capm}{\morstuff{cap}}
\newcommand{\cupm}{\morstuff{cup}}

\newcommand{\capp}{\morstuff{pcap}}
\newcommand{\cupp}{\morstuff{pcup}}

\newcommand{\tcapm}{\morstuff{tup}}
\newcommand{\tcupm}{\morstuff{tdown}}

\newcommand{\tcapp}{\morstuff{tup}}
\newcommand{\tcupp}{\morstuff{tdown}}

\newcommand{\slnbasis}[2]{\setstuff{C}_{#1}^{#2}}
\newcommand{\glnbasis}[2]{\setstuff{CP}_{#1}^{#2}}
\newcommand{\sonbasis}[2]{\setstuff{P}_{#1}^{#2}}

\newcommand{\symweb}[1][{\sln[2]}]{\twocatstuff{S}\twocatstuff{W}\catstuff{eb}(#1)}

\usepackage[all]{xy}
\usepackage{tikz}
\usetikzlibrary{cd}
\usetikzlibrary{decorations}
\usetikzlibrary{decorations.markings}
\usetikzlibrary{decorations.pathreplacing}
\usetikzlibrary{decorations.pathmorphing}
\usetikzlibrary{arrows.meta,shapes,positioning,matrix,calc}
\usetikzlibrary{shapes.callouts}
\tikzset{anchorbase/.style={baseline={([yshift=-0.5ex]current bounding box.center)}},
tinynodes/.style={font=\tiny,text height=0.25ex,text depth=0.05ex},
smallnodes/.style={font=\scriptsize,text height=0.75ex,text depth=0.15ex},
crossline/.style={preaction={draw=white,line width=10.0pt,-},preaction={draw=black,line width=1.8pt,-}},
usual/.style={line width=2.0,color=black},
symedge/.style={line width=2.0,double,color=tomato},
phantom/.style={line width=2.0,color=magenta,densely dashed},
}
\tikzstyle directed=[postaction={decorate,decoration={markings,mark=at position #1 with {\arrow[line width=0.5mm, black]{>}}}}]
\tikzstyle rdirected=[postaction={decorate,decoration={markings,mark=at position #1 with {\arrow[line width=0.5mm, black]{<}}}}]

\def\NewTheorem#1{%
\newaliascnt{#1}{equation}%
\newtheorem{#1}[#1]{#1}%
\aliascntresetthe{#1}%
\expandafter\def\csname #1autorefname\endcsname{#1}%
}
\def\equationautorefname~#1\null{(#1)\null}

\numberwithin{equation}{subsection}

\NewTheorem{Proposition}
\NewTheorem{Theorem}
\NewTheorem{Corollary}
\AtEndEnvironment{Corollary}{\null\hfill$\square$}%
\NewTheorem{Lemma}
\theoremstyle{definition}
\NewTheorem{Definition}
\AtEndEnvironment{Definition}{\null\hfill$\Diamond$}%
\NewTheorem{Classification Problem}
\AtEndEnvironment{Classification Problem}{\null\hfill$\Diamond$}%
\NewTheorem{Notation}
\AtEndEnvironment{Notation}{\null\hfill$\Diamond$}%
\NewTheorem{Example}
\AtEndEnvironment{Example}{\null\hfill$\Diamond$}%
\NewTheorem{Examples}
\AtEndEnvironment{Examples}{\vskip-10mm\null\hfill$\Diamond$}%

\theoremstyle{remark}
\NewTheorem{Remark}
\AtEndEnvironment{Remark}{\null\hfill$\Diamond$}%
\NewTheorem{Assumption}

\usepackage[hypertexnames=false]{hyperref}
\usepackage{bookmark}
\hypersetup{
pdftoolbar=true,
pdfmenubar=true,
pdffitwindow=false,
pdfstartview={FitH},
pdftitle={On rank one 2-representations of web categories},
pdfauthor={Daniel Tubbenhauer},
pdfsubject={},
pdfcreator={Daniel Tubbenhauer},
pdfproducer={Daniel Tubbenhauer},
pdfkeywords={},
pdfnewwindow=true,
colorlinks=true,
linkcolor=mydarkblue,
citecolor=teal,
filecolor=magenta,
urlcolor=orchid,
linkbordercolor=lava,
citebordercolor=teal,
urlbordercolor=orchid,
linktocpage=true
}

\newcommand{\nnfootnote}[1]{%
\begin{NoHyper}
\renewcommand\thefootnote{}\footnote{#1}%
\addtocounter{footnote}{-1}%
\end{NoHyper}
}

\def\makeautorefname#1#2{\csdef{#1autorefname}{#2}}
\makeautorefname{section}{Section}%
\makeautorefname{subsection}{Section}%
\makeautorefname{subsubsection}{Section}%

\begin{document}
\title[On rank one 2-representations of web categories]{On rank one 2-representations of web categories}
\author[D. Tubbenhauer]{Daniel Tubbenhauer}

\address{D.T.: The University of Sydney, School of Mathematics and Statistics F07, Office Carslaw 827, NSW 2006, Australia, \href{http://www.dtubbenhauer.com}{www.dtubbenhauer.com}, https://orcid.org/0000-0001-7265-5047}
\email{daniel.tubbenhauer@sydney.edu.au}

\begin{abstract}
We classify rank one 2-representations of $\mathrm{SL}_{2}$, 
$\mathrm{GL}_{2}$ and 
$\mathrm{SO}_{3}$ web categories.
The classification is inspired by similar results about quantum groups, given by reducing the problem to the classification of 
bilinear and trilinear forms, and is formulated such that
it can be adapted to other web categories.
\end{abstract}

\nnfootnote{\textit{Mathematics Subject Classification 2020.} Primary: 18M05, 18M30; Secondary: 22E46, 22E47.}
\nnfootnote{\textit{Keywords.} Webs, 2-representations, bilinear and trilinear forms.}

\addtocontents{toc}{\protect\setcounter{tocdepth}{1}}

\maketitle

\tableofcontents


\section{Introduction}\label{S:Introduction}


We \changed{give a} classification of simple transitive
2-representations of web categories on $\C$-vector spaces.
This classification \ochanged{builds upon and extends} results in the classification of quantum groups and Hopf algebras. 
The main point is that, even in this very restricted and semisimple setting,
there are infinitely many such 2-representations 
and \ochanged{these} are essentially impossible to classify explicitly.


\subsection{Motivation and results}\label{SS:Introduction}


\emph{Classification} is a central topic in all of mathematics. In representation theory the most important classification problem is to construct and compare all simple representations. In higher representation theory, 
an offspring of \emph{categorification} that originates in seminal papers such as \cite{ChRo-sl2-categorification}, \cite{EtNiOs-fusion-categories}, \cite{KhLa-cat-quantum-sln-first}, \cite{Os-module-categories} or \cite{MaMi-cell-2-reps}, the most crucial classification problem is about the appropriate analog 
of simple representations. For example, given a favorite monoidal category, one can ask whether one can classify its simplest possible module categories.
The favorite categories of our choice in this note are certain \emph{diagram categories}, simplest possible will mean \emph{simple transitive} and classification will mean \emph{reduction} of the original problem to linear algebra.

Note\changed{, however,} that linear algebra can still be \ochanged{arbitrarily} complicated.
The problem of classifying symmetric and alternating \emph{bilinear forms} is well-known and has a very pleasant answer. Less well-known 
but still doable and nice is the classification of all bilinear forms.
On the other hand, the classification of \emph{trilinear forms}
is tractable for small dimensions only, even if one restricts to
symmetric or alternating forms: the classification problem is ``wilder than wild'' \cite{BeSe-matrix-problems}. However, for small dimensions there is indeed a classification of trilinear forms, see e.g. \cite{ThCh-trilinear}, \cite{Ch-invariants-trilinear} or \cite{Th-trilinear}.

In this note we will see a similar behavior for the following
\emph{web categories}: the category of $\sln=\sln(\C)$ webs $\web$, 
the category of $\gln=\gln(\C)$ webs $\web[{\gln}]$ and the 
category of $\son=\son(\C)$ webs $\web[{\son}]$, and quantum versions \changed{for which the $\qpar$ in the notation will appear}.
(That we discuss $\son$ webs and not the very similar $\sln[3]$ webs 
has historical reasons, see \autoref{R:RepsSO3HistRemark} below.)
The classification problem we have in mind for these categories
is to study the easiest form of actions of these categories on $\vect$, the category of finite dimensional $\C$-vector spaces. In the language of \cite{MaMi-transitive-2-reps}, we want to classify \emph{rank one simple transitive 2-representations} of these web categories. \changed{(Let us point out that rank one simple transitive 2-representations are not semisimple in general, but for web categories they are by \autoref{L:RepsSemisimple}.)}

For all of these web categories we give a classification of such 
2-representations. 
\changed{The classification takes a certain form as outlined in \autoref{C:RepsClassification}. Roughly,} 
we construct 2-representations from \ochanged{linear algebra inputs} such as matrices 
and tensors. Second, the equivalence classes of these 2-representations 
are given by an explicit linear algebra condition on matrices and tensors such as congruence.
Finally, we argue that any rank one simple transitive 2-representation
is of the form constructed in the first step.

But how explicit our classification is varies drastically:
\begin{enumerate}[label=$\triangleright$]

\item For $\web$ the classification is similar to the classification of 
bilinear forms and has therefore a short-and-sweet answer, see 
\autoref{T:RepsSL2Main}.

\item For $\web[{\gln}]$ trilinear forms make their appearance.
However, as we will see, the appearing trilinear forms are on 
$\C$-vector spaces of small dimensions so we still get a good answer, see
\autoref{T:RepsGL2Main}.

\item For $\web[{\son}]$ honest trilinear forms appear which makes us 
believe that there is no good (this could e.g. mean listable) answer, see \autoref{T:RepsSO3Main}.

\end{enumerate}
These three web categories are 
semisimple but have infinitely many isomorphism classes
of simple objects.
As we will see, in all cases there are infinitely many equivalence
classes of rank one simple transitive 2-representations. This is 
very different from the situation of semisimple categories with finitely
many simple objects where some form of Ocneanu rigidity ensures that there
are only finitely many simple transitive 2-representations.

In \autoref{P:ComplexitySL2} we also show that the classification (of rank one simple transitive 2-representations) 
for $\web$ (and $\web[{\gln}]$) implies the classification of bilinear forms, and in fact, the classification is a tame problem.
For $\web[{\son}]$ we are not able to determine the precise characterization of the complexity of the classification problem. However, for a modification of $\web[{\son}]$ we show that the classification implies the classification of trilinear forms, see \autoref{P:ComplexitySO3}. In fact, the classification problem 
for the variant of $\web[{\son}]$ is strictly more difficult than any wild problem in classical representation theory, see \autoref{T:ComplexityDifficult}.


\subsection{A few extra comments}\label{SS:IntroductionComments}


We finish the introduction with a few \ochanged{remarks}.

\begin{Remark}\label{R:IntroWebs}
All web categories in this paper are monoidally equivalent to
representation theoretical categories. However, since one of our main points is
to use diagrammatics, we think of these 
as web categories instead of their representation theoretical counterparts.

Along the same lines, we \ochanged{would} like to point out that similar
results have been obtained in other fields although the translation 
is not completely straightforward. The connection was in fact the starting point 
for this note. The methods presented in these papers are different 
\ochanged{from} what we do in this note; in particular, we take the diagrammatic approach
and make the classification results more explicit, see e.g. 
\autoref{L:RepsSL2List}.

For example, see \cite{Bi-repcat-qgroup} or \cite{EtOs-module-categories-slq2} for $\sln[2]$, \cite{Mr-quantum-gl2} 
for $\gln$, and \cite{Mr-quantum-so3} or \cite{EvPu-so3-classification} for $\son$. See also \cite{Oh-sl3-cat-reps} and \cite{NeYa-quantum-sun} for the $\sln[k]$ family.
\end{Remark}

\begin{Remark}\label{R:IntroOtherWebs}
In this remark \emph{complexity} is meant with respect to 
classification of rank one 
simple transitive 2-representations\changed{, and we use it as an informal estimate of difficulty. We give some details later in \autoref{S:Complexity}.}
\begin{enumerate}

\item Consider the following list:
\begin{gather*}
\begin{gathered}
\sln[2],\gln
\\[-0.1cm]
\begin{tabular}{C||C|C|C|C|C|C|C}
n & 0 & 1 & 2 & 3 & 4 & 5 & 6 \\
\hline
\dim_{\C} & 1 & 0 & 1 & 0 & 2 & 0 & 5 \\
\end{tabular}
\end{gathered}
,\quad
\begin{gathered}
\sln[3]
\\[-0.1cm]
\begin{tabular}{C||C|C|C|C|C|C|C}
n & 0 & 1 & 2 & 3 & 4 & 5 & 6 \\
\hline
\dim_{\C} & 1 & 0 & 1 & 1 & 2 & 3 & 6 \\
\end{tabular}
\end{gathered}
,\\
\begin{gathered}
\son
\\[-0.1cm]
\begin{tabular}{C||C|C|C|C|C|C|C}
n & 0 & 1 & 2 & 3 & 4 & 5 & 6 \\
\hline
\dim_{\C} & 1 & 0 & 1 & 1 & 3 & 6 & 15 \\
\end{tabular}
.
\end{gathered}
\end{gather*}
These lists are the maximal appearing dimensions $b_{n}^{\ast}$ of the hom-spaces in $\sln[2]$, $\sln[3]$ and $\son$ webs, respectively, for webs with $n$ boundary points to the empty web. We have $b_{n}^{\sln}\leq b_{n}^{\sln[3]}\leq b_{n}^{\son}$, but 
on the other hand \cite[Theorem 1.4.(a)]{CoOsTu-growth} gives
\begin{gather*}
\lim_{n\to\infty}\sqrt[n]{b_{n}^{\sln}}=2
\,\colorbox{tomato!50}{$<$}\,
3
=\lim_{n\to\infty}\sqrt[n]{b_{n}^{\sln[3]}}=\lim_{n\to\infty}\sqrt[n]{b_{n}^{\son}}.
\end{gather*}
This justifies the complexity jump from $\sln[2]$ 
to $\son$ webs, and probably 
indicates that $\sln[3]$ and $\son$ webs are of the same complexity.
(\cite{Bi-asymptotic-lie} gives more precise formulas for the asymptotics of the numbers $b_{n}^{\ast}$ but we do not need them here.)

Note also that $b_{n}^{\sln}=b_{n}^{\gln}$ so their complexity is roughly the 
same which indeed matches what we will see in \autoref{T:RepsSL2Main} and \autoref{T:RepsGL2Main}.

\item In general we expect the complexity 
of $\sln[n]$ (or $\gln[n]$) webs, as 
in e.g. \cite{MuOhYa-webs} or \cite{CaKaMo-webs-skew-howe}, to be 
equal to \ochanged{or} higher than for $\son$, so likely unsolvable in 
a precise sense. However, as pointed out in \cite{Th-trilinear}, 
$(4,4,6)$ trilinear forms are classifiable and $4$, $4$, $6$ 
are the dimensions of the nontrivial fundamental $\sln[4]$-representations.
Thus, there might be something \ochanged{that can be said} for $\sln[4]$ webs similar to 
what we do at the end of \autoref{SS:RepsSO3Proof}.

\item All categories in this note are semisimple. A good question 
is to address the nonsemisimple case where one could expect 
\emph{cell theory} as in e.g. 
\cite{Gr-structure-semigroups}, \cite{MaMi-cell-2-reps} or 
\cite{Tu-sandwich-cellular} to play a role. As usual one should 
expect a nontrivial complexity jump from the semisimple 
to the nonsemisimple cases.

\end{enumerate}
(In (b) we write \emph{$(p,q,r)$ trilinear form} for a trilinear form 
on $\C$-vector spaces of these dimensions.)
\end{Remark}

\begin{Remark}\label{R:IntroColors}
Two technical remarks:
\begin{enumerate}

\item Some calculations were done with machine help and the \ochanged{code used} 
can be found on \cite{Tu-code-2-reps-webs}. 
\ochanged{Any necessary erratum can be found on the same website.}

\item This note is readable in black-and-white since
the colors we use are a visual aid only.

\end{enumerate}
(Please email me if you find mistakes or have comments; I aim to then update \cite{Tu-code-2-reps-webs}.)
\end{Remark}

\noindent\textbf{Acknowledgments.} D.T. was 
sponsored by the ARC Future Fellowship FT230100489. Although D.T. gets 
more pathetic by the second, they \ochanged{would} like to thank Victor Ostrik for 
pointing them to the paper by Bichon which was the starting 
point of this project. They also thank Joanna Meinel
for nice walks \ochanged{during which many good questions were asked}, 
Anne Dranowski and Pedro Vaz for useful discussions,
and Noah Snyder for helping \ochanged{to clarify the historical facts}. D.T. also acknowledges ChatGPT for its assistance with proofreading. Special thanks to the referee\tchanged{s} for a careful reading of the manuscript and very helpful comments and suggestions.


\section{Background on 2-representations}\label{S:Reps}


We will briefly recall some notions from 2-representation theory. \changed{Below, and throughout, we usually count objects up to isomorphism but drop the `up to isomorphism' for brevity.} 

\begin{Remark}\label{D:RepsReformulates}
Nothing in this section is new, but reformulated
compared to our main 
sources \cite{EtGeNiOs-tensor-categories},  
\cite{MaMaMiTuZh-bireps} and \cite{MaMaMiTuZh-soergel-2reps}.
Our reformulation stems from that we only consider the semisimple 
case, which is just a special case of what the above literature discusses.
For categories with countably many 
simple object see also \cite{Ma-infinite-finitary} for a helpful reference.

The reader is also \ochanged{referred} to \cite{MaLa-categories} or \cite{EtGeNiOs-tensor-categories} 
for standard notions such as monoidal categories.
\end{Remark}

Let $\K$ be \changed{a field}. Everything in this section will be over $\K$.
As basic notation we use:
\begin{enumerate}[label=$\triangleright$]

\item Categories are denoted using bold font \changed{such as $\catstuff{V}$ or $\catstuff{W}$}, 
while monoidal categories ($\leftrightsquigarrow$ 2-categories with one object) 
are denoted by e.g. \changed{$\twocatstuff{C}$ or $\twocatstuff{D}$}.

\item Objects are denoted by \changed{$\obstuff{X}$, $\obstuff{Y}$} 
etc. and morphisms by \changed{$\morstuff{f}$, $\morstuff{g}$} etc.

\item Categorical composition is denoted by $\vcirc$ and monoidal composition 
by $\hcirc$. The monoidal unit is denoted by $\munit$, and identity morphism 
as $\morstuff{id}_{\obstuff{X}}$ or simplified as 
$\morstuff{id}$. The following picture summarizes our diagrammatic composition rules:
\begin{gather*}
(\idmor\hcirc\morstuff{g})\vcirc(\morstuff{f}\hcirc\idmor) 
=
\begin{tikzpicture}[anchorbase,scale=0.22,tinynodes]
\draw[spinach!35,fill=spinach!35] (-5.5,2) rectangle (-2,0.5);
\draw[spinach!35,fill=spinach!35] (-0.5,3) rectangle (3,4.5);
\draw[thick,densely dotted] (-5.5,2.5) node[left,yshift=-2pt]{$\vcirc$} 
to (3.5,2.5) node[right,yshift=-2pt]{$\vcirc$};
\draw[thick,densely dotted] (-1.25,5) 
node[above,yshift=-2pt]{$\hcirc$} to (-1.25,0) node[below]{$\hcirc$};
\draw[usual] (-5,2) to (-5,3.75) node[right]{$\dots$} to (-5,5);
\draw[usual] (-2.5,2) to (-2.5,5);
\draw[usual] (0,0) to (0,1.25) node[right]{$\dots$} to (0,3);
\draw[usual] (2.5,0) to (2.5,3);
\draw[usual] (-5,0) to (-5,0.25) node[right]{$\dots$} to (-5,0.5);
\draw[usual] (-2.5,0) to (-2.5,0.5);
\draw[usual] (2.5,4.5) to (2.5,5);
\draw[usual] (0,4.5) to (0,4.75) node[right]{$\dots$} to (0,5);
\node at (-3.75,1.0) {$\morstuff{f}$};
\node at (1.25,3.5) {$\morstuff{g}$};
\end{tikzpicture}
=
\begin{tikzpicture}[anchorbase,scale=0.22,tinynodes]
\draw[spinach!35,fill=spinach!35] (-5.5,1.75) rectangle (-2,3.25);
\draw[spinach!35,fill=spinach!35] (-0.5,1.75) rectangle (3,3.25);
\draw[usual] (-5,3.25) to (-5,4.125) node[right]{$\dots$} to (-5,5);
\draw[usual] (-2.5,3.25) to (-2.5,5);
\draw[usual] (0,0) to (0,0.625) node[right]{$\dots$} to (0,1.75);
\draw[usual] (2.5,0) to (2.5,1.75);
\draw[usual] (-5,0) to (-5,0.625) node[right]{$\dots$} to (-5,1.75);
\draw[usual] (-2.5,0) to (-2.5,1.75);
\draw[usual] (2.5,3.25) to (2.5,5);
\draw[usual] (0,3.25) to (0,4.125) node[right]{$\dots$} to (0,5);
\node at (-3.75,2.25) {$\morstuff{f}$};
\node at (1.25,2.25) {$\morstuff{g}$};
\end{tikzpicture}
=
\\
\begin{tikzpicture}[anchorbase,scale=0.22,tinynodes]
\draw[spinach!35,fill=spinach!35] (5.5,2) rectangle (2,0.5);
\draw[spinach!35,fill=spinach!35] (0.5,3) rectangle (-3,4.5);
\draw[thick, densely dotted] (5.5,2.5) 
node[right,yshift=-2pt]{$\vcirc$} to (-3.5,2.5) node[left,yshift=-2pt]{$\vcirc$};
\draw[thick, densely dotted] (1.25,5) 
node[above,yshift=-2pt]{$\hcirc$} to (1.25,0) node[below]{$\hcirc$};
\draw[usual] (2.5,2) to (2.5,3.75) node[right]{$\dots$} to (2.5,5);
\draw[usual] (5,2) to (5,5);
\draw[usual] (0,0) to (0,3);
\draw[usual] (-2.5,0) to (-2.5,1.25) node[right]{$\dots$} to (-2.5,3);
\draw[usual] (2.5,0) to (2.5,0.25) node[right]{$\dots$} to (2.5,0.5);
\draw[usual] (5,0) to (5,0.5);
\draw[usual] (-2.5,4.5) to (-2.5,4.75) node[right]{$\dots$} to (-2.5,5);
\draw[usual] (0,4.5) to (0,5);
\node at (3.75,1.0) {$\morstuff{g}$};
\node at (-1.25,3.5) {$\morstuff{f}$};
\end{tikzpicture}
=(\morstuff{f}\hcirc\idmor)\vcirc(\idmor\hcirc\morstuff{g}).
\end{gather*}

\item For functors calligraphic font is used, for example \changed{$\functorstuff{M}$ or $\functorstuff{N}$}.

\end{enumerate}

\begin{Remark}\label{R:RepsStrict}
Before we get started, let us point out that \ochanged{all} our \changed{higher} categories 
and structures are strict. This turns out to be no restriction for us
and the usual strictification theorems apply, see for example 
\cite[Section 1.4]{GoPoSt-coherence-tricategories} or \cite[Section 4.2]{Po-general-coherence} for quite general strictification theorems.

\changed{In contrast, the notion of equivalence (and
generally of morphism) for 2-representations is not strict. See also \cite{Sn-strict} for examples why one has to allow nonstrict equivalences.}
\end{Remark}

\begin{Definition}\label{D:RepsSemisimpleDef}
A nonzero object of a $\K$-linear category is \emph{simple} if \ochanged{all its subobjects are either zero objects or isomorphisms}.

By a \emph{semisimple category} we mean a 
\changed{$\K$-linear additive category} such that (i) every object decomposes as a finite direct sum of simple objects, and (ii) if $\morstuff{f}\colon\obstuff{S}_{1}\to\obstuff{S}_{2}$ 
and $\morstuff{g}\colon\obstuff{S}_{2}\to\obstuff{S}_{3}$ are nonzero 
morphisms between simple objects $\obstuff{S}_{i}$, \ochanged{then $\morstuff{g}\vcirc\morstuff{f}$ is also nonzero.}	
\end{Definition}

A semisimple category \ochanged{may} have 
infinitely many simple objects.
The following will be used silently.

\begin{Lemma}\label{L:RepsSemisimpleDef}
We have:
\begin{enumerate}

\item \changed{A semisimple category is idempotent complete.}

\item A semisimple category is $\K$-linear finite length abelian.

\item $\K$-linear functors between $\K$-linear additive categories are additive.

\end{enumerate}
\end{Lemma}

\begin{proof}
\changed{Well-known, so we will be brief. First, (c) is an exercise, while (b) is proven in, for example, \cite[Theorem C.6]{TuVi-monoidal-tqft}. (a) follows from the following argument (due to the referee):
For (a) it suffices to show that an idempotent in the endomorphism algebra of a
simple object splits. And indeed, 
if $\morstuff{e}$ is a nontrivial (neither zero nor the identity) idempotent in the endomorphism algebra of a simple object $\obstuff{X}$, then so is $\idmor_{\obstuff{X}}-\morstuff{e}$. However, the composite of these two is zero, violating (ii) of \autoref{D:RepsSemisimpleDef}.}
\end{proof}

\begin{Remark}\label{R:RepsSame}
By \autoref{L:RepsSemisimpleDef}.(a)\changed{, in the semisimple case,} the theories 
presented in \cite{EtGeNiOs-tensor-categories} on the one hand, and 
\cite{MaMaMiTuZh-bireps} and \cite{MaMaMiTuZh-soergel-2reps} on the other hand are essentially the same (there are subtle differences but they do not play any role for us).
\end{Remark}

\begin{Definition}\label{D:RepsFusion}
\changed{An \emph{essentially fusion category $\twocatstuff{C}$}} is a 
semisimple \ochanged{rigid monoidal} category with countably \ochanged{many simple objects} 
and finite dimensional morphism spaces.
\end{Definition}

\begin{Example}\label{E:RepsFusion}
The category $\vect[\K]$ of finite dimensional $\K$-vector spaces 
is a prototypical example of \changed{an essentially fusion} category. A more exciting 
example is $\rep\big(\matsln\big)$, complex finite dimensional $\matsln$-representations, and its relatives that we will 
discuss in the sections below. The simple objects, up to equivalence, 
in the category 
$\rep\big(\matsln\big)$ are $\{\obstuff{S}_{k}=\mathrm{Sym}^{k}(\C^{2})|k\in\Z_{\geq 0}\}$.
\end{Example}

\begin{Remark}\label{D:RepsSemisimple}
\changed{In the language of \cite{MaMaMiTuZh-soergel-2reps}
\autoref{D:RepsFusion} translates to 
what is called locally semisimple quasi-fiat one object 2-category 
in that paper with one difference: \autoref{D:RepsFusion} allows countably (finite or infinite) many simple objects while locally semisimple quasi-fiat one object 2-categories always have finitely many simple objects.}
\end{Remark}

\begin{Lemma}\label{L:RepsKrullSchmidt}
\changed{An essentially fusion} category $\twocatstuff{C}$ is Krull--Schmidt.
\end{Lemma}

\begin{proof}
A finite length abelian category is Krull--Schmidt, so 
\ochanged{\autoref{L:RepsSemisimpleDef}.(b) proves the claim}.
\end{proof}

Let $\fincat=\fincat_{\K}^{f}$ denote the 2-category of finitary categories, $\K$-linear functors and natural \changed{transformations}, see 
\cite[Definition 2.12]{MaMaMiTuZh-bireps}. For the purpose of this 
paper it is enough to know that $\vect[\K]\in\fincat$.
\changed{By $\vcirc$-ideal we mean an ideal with respect to the operation $\vcirc$, while $\vcirc$-$\hcirc$-ideal is meant with respect to both operations $\vcirc$ and $\hcirc$ separately.}

\begin{Definition}\label{D:RepsReps}
Let $\twocatstuff{C}$ be as in \autoref{D:RepsFusion}.
\begin{enumerate}

\item A \emph{(finitary) 2-representation} $\trep$ of $\twocatstuff{C}$ is a 
$\K$-linear monoidal functor $\trep\colon\twocatstuff{C}\to\End_{\ffincat}(\catstuff{V})$ for 
$\catstuff{V}\in\fincat$.

\item The \emph{rank}, denoted by $\rank\trep$, of such a functor $\trep$ is the number of indecomposable objects in $\catstuff{V}$.	

\item We call a 2-representation \emph{semisimple} if 
the target category $\catstuff{V}\in\fincat$ is semisimple.

\item Such a functor $\trep$ is called \emph{simple transitive} if it has no proper $\twocatstuff{C}$-stable $\vcirc$-ideals, meaning that every $\vcirc$-ideal $I\subset\catstuff{V}$ with $\trep(\obstuff{X})(I)\subset I$ for all $\obstuff{X}\in\twocatstuff{C}$ 
is either zero or $\catstuff{V}$.

\end{enumerate}
(The 2-representation in (d) are simple, by the definition above, and transitive by \autoref{L:RepsTransitive} below, hence the name.)
\end{Definition}

\begin{Remark}\label{R:RepsReps}
There is also the notion of a 
\emph{(finitary) module category}. Similarly as for representations 
and modules, these notions are equivalent in the appropriate sense. 
We leave it to the reader to spell out the definitions, and we use them \changed{interchangeably}. \tchanged{For example, in \autoref{D:RepsEquivalence} below the horizontal arrows are to be read as module category notation.}
\end{Remark}

\begin{Definition}\label{D:RepsEquivalence}
Two 2-representations $\trep\colon\twocatstuff{C}\to\End_{\ffincat}(\catstuff{V})$ 
and $\trep[N]\colon\twocatstuff{C}\to\End_{\ffincat}(\catstuff{W})$ 
are \emph{equivalent}, written $\trep\equirep\trep[N]$, if there is an
equivalence $\functorstuff{F}\colon\catstuff{V}\to\catstuff{W}\in\fincat$ 
such that
\begin{gather*}
\begin{tikzpicture}[anchorbase,->,>=stealth',shorten >=1pt,auto,node distance=3cm,
thick,main node/.style={font=\sffamily\Large\bfseries}]
\node[main node] (1) {$\twocatstuff{C}$};
\node[main node] (2) [right of=1,xshift=-1.0cm] {$\catstuff{V}$};
\node[main node] (3) [below of=1,yshift=2.0cm] {$\twocatstuff{C}$};
\node[main node] (4) [right of=3,xshift=-1.0cm] {$\catstuff{W}$};
\path[every node/.style={font=\sffamily\small}]
(1) edge node[above]{$\trep$} (2)
(3) edge node[below]{$\trep[N]$} (4)
(1) edge[double equal sign distance,-] (3)
(2) edge node[right]{$\functorstuff{F}$} (4);
\end{tikzpicture}
\end{gather*}
is a commutative diagram \changed{up to a coherent natural isomorphism $\trep[N](\obstuff{X})\big(\functorstuff{F}(\obstuff{V})\big)\xRightarrow{\cong}\functorstuff{F}\big(\trep[M](\obstuff{X})(\obstuff{V})\big)$}.
\end{Definition}

\begin{Lemma}\label{L:RepsTransitive}
Any simple transitive 2-representation is \emph{transitive}, meaning 
generated\changed{, taking direct sums and direct summands,} by the action of $\trep$.
\end{Lemma}

\begin{proof}
This is \cite[Lemma 4]{MaMi-transitive-2-reps}.
\end{proof}

\begin{Example}\label{E:RepsFusionSimple}
The category $\rep\big(\matsln\big)$ acts on itself by tensoring, that is, $\trep(\obstuff{S}_{k})$ is the endofunctor of (left) tensoring with $\obstuff{S}_{k}$. 
The only $\vcirc$-$\hcirc$-ideals in $\rep\big(\matsln\big)$ are 
zero or the category itself. This follows since
\changed{$\rep\big(\matsln\big)$ is semisimple and we have}
$\C^{2}\hcirc\obstuff{S}_{k}\cong\obstuff{S}_{k+1}\oplus\obstuff{S}_{k-1}$ 
for $k\in\Z_{\geq 1}$. Hence, $\rep\big(\matsln\big)$ is simple transitive.
Thus, since $\rank\rep\big(\matsln\big)=\infty$, $\rep\big(\matsln\big)$ is an infinite rank simple transitive 2-representation of itself.
\end{Example}

Note that 2-representations of \changed{essentially fusion} categories are in general not 
semisimple. As an example consider $\vect[\K]$ which can act on 
any $\K$-linear abelian category. \changed{(This action is unique up to the equivalence in \autoref{D:RepsEquivalence}.)} Hence, the following \ochanged{result} is remarkable and key for this paper:

\begin{Lemma}\label{L:RepsSemisimple}
\tchanged{Any simple transitive 2-representation $\trep$ of \changed{an essentially fusion} category $\twocatstuff{C}$ with $\rank\trep<\infty$ is semisimple, meaning that $\catstuff{V}$ is semisimple.}
\end{Lemma}

\begin{proof}
A direct adaption of \cite[Proposition 2.16]{MaMaMiTuZh-soergel-2reps}.
\end{proof}

\begin{Lemma}\label{L:RepsSemisimpleRankOne}
Any rank one simple transitive 2-representation of \changed{an essentially fusion} category $\twocatstuff{C}$ is on $\vect[\K]$.
\end{Lemma}

\begin{proof}
Since $\vect[\K]$ is the only semisimple category with one simple object, this 
follows from \autoref{L:RepsSemisimple}.
\end{proof}

Let $\twocatstuff{C}$ be \changed{an essentially fusion} category.
Recall that a \emph{fiber functor} $\mrep\colon\twocatstuff{C}\to\vect[\K]$ 
is an exact faithful monoidal functor. \changed{We write $\mrep(\obstuff{X}|\obstuff{X}\in\twocatstuff{C})$ for the full subcategory whose objects are direct sums of direct summands of
objects of the form $\mrep(\obstuff{X})$, for 
$\obstuff{X}\in\twocatstuff{C}$.} The following, very easy, lemma is another key fact:

\begin{Lemma}\label{L:RepsFiber}
We have the following.
\begin{enumerate}

\item Any fiber functor $\mrep\colon\twocatstuff{C}\to\vect[\K]$ of \changed{an essentially fusion} category $\twocatstuff{C}$ give\changed{s} rise to a semisimple rank one 2-representation $\trep$.

\item A fiber functor functor $\mrep\colon\twocatstuff{C}\to\vect[\K]$ of \changed{an essentially fusion} category $\twocatstuff{C}$ give\changed{s} rise to a simple transitive 2-representation $\trep$ if and only if \changed{$\mrep(\obstuff{X}|\obstuff{X}\in\twocatstuff{C})$} does not have any nontrivial $\vcirc$-$\hcirc$-ideals.

\end{enumerate}

\end{Lemma}

\begin{proof}
\textit{(a).}
For $\obstuff{X}\in\twocatstuff{C}$ set $\trep(\obstuff{X})$ to be the 
endofunctor of tensoring with the $\K$-vector space $\mrep(\obstuff{X})$. 
One \ochanged{can directly verify} that this defines a 2-representation.

\textit{(b).} If $\trep$ is simple transitive, then there cannot be any 
nontrivial $\vcirc$-$\hcirc$-ideals by the construction of $\trep$ in (a).
Conversely, if there are no nontrivial $\vcirc$-$\hcirc$-ideals then semisimplicity, that is \autoref{L:RepsSemisimple}, \tchanged{implies} that $\trep$ is  
simple transitive.
\end{proof}

With \autoref{L:RepsFiber}.(a) in mind, 
we also say \emph{fiber 2-representation} instead of fiber functor. These 
are always of rank one, by definition, but the converse might be false. 
(A rank one 2-representation has no reason to be faithful in general.)
\changed{If}
the condition \autoref{L:RepsFiber}.(b) is satisfied for a fiber functor 
$\mrep\colon\twocatstuff{C}\to\vect[\K]$, then we call $\mrep$ 
a \emph{simple transitive fiber 2-representation}.

\begin{Example}\label{E:RepsFiber2Rep}
The action \tchanged{from} \autoref{E:RepsFusionSimple} 
is not a fiber 2-representation. But composition with the 
forgetful functor $\rep\big(\matsln\big)\to\vect$ 
defines a (simple transitive) fiber 2-representation.
\end{Example}

\begin{Classification Problem}\label{C:RepsClassification}
The classification of simple transitive 2-representa\-tions of a given 
$\twocatstuff{C}$ is one of the main problems of the theory, and akin 
to classifying simple representations of groups or algebras. This is 
justified by the categorical analog of the Jordan--H{\"o}lder theorem, 
see \cite[Section 3.5]{MaMi-transitive-2-reps}.

For us such a classification is optimally given by:
\begin{enumerate}[label=$\triangleright$]

\item The existence of certain explicitly constructed simple transitive 2-representa\-tions. \textbf{(Existence)}

\item The comparison of these with a computable condition. \textbf{(Non-redundant)}

\item A proof that all simple transitive 2-representations are of the particular form. \textbf{(Complete)}

\end{enumerate}
In this paper we restrict to the subproblem of classifying 
simple transitive rank one (or fiber) 2-representations. As we will see, 
even this subproblem can get \changed{arbitrarily} difficult, and we will sometimes 
only give part of the list above.
\end{Classification Problem}

\begin{Remark}\label{R:RepsClassification}
\autoref{C:RepsClassification} is not meant as a definition.
\end{Remark}

\begin{Example}\label{E:RepsClassification}
Keeping \autoref{R:RepsSame} in mind, the paper \cite{EtOs-module-categories-slq2} classifies simple transitive 2-representations 
of $\rep\big(\matsln\big)$ of finite rank. The classification is quite difficult, and we will discuss the 
much simpler classification 
of simple transitive fiber 2-representations $\rep\big(\matsln\big)\to\vect[\C]$ 
in \autoref{S:RepsSL2}. It turns out that in this case 
all rank one simple transitive 2-representations come from fiber 
functors.
\end{Example}

\begin{Lemma}\label{L:RepsSemisimpleOcneanu}
\changed{An essentially fusion} category $\twocatstuff{C}$ with finitely many simple 
objects has only finitely many simple transitive 2-representations up to $\equirep$.
\end{Lemma}

\begin{proof}
\changed{We point out that an essentially fusion category with finitely many simple objects is a fusion category in the usual sense as, for example, in \cite[Chapter 9]{EtGeNiOs-tensor-categories}. Then the claim follows from} Ocneanu rigidity as e.g. 
in \cite[Proposition 3.4.6 and Corollary 9.1.6]{EtGeNiOs-tensor-categories}.
\end{proof}

With contrast to \autoref{L:RepsSemisimpleOcneanu} we have:

\begin{Theorem}\label{T:RepsSemisimpleOcneanu}
\changed{An essentially fusion} category $\twocatstuff{C}$ can have infinitely many nonequivalent simple transitive rank one 2-representations.
\end{Theorem}

\begin{proof}
\changed{By the examples discussed in the next sections; see, for example, \autoref{T:RepsSL2Main}.}
\end{proof}	

For $\vect[\K]$ we use its 
unique braiding (= the flip map).
To finish this section, and relevant for our examples \ochanged{(we say braided and symmetric instead of braided monoidal and symmetric monoidal)}:

\begin{Definition}\label{D:RepsBraidedReps}
Assume that the acting category $\twocatstuff{C}$ is braided.
A fiber 2-representation is \emph{braided} if it is given by a braided functor.
\end{Definition}

One should not expect fiber 2-representations to have interesting braidings:

\begin{Lemma}\label{L:RepsBraidedReps}
A braided fiber 2-representation is given by a symmetric functor.
\end{Lemma}

\begin{proof}
\changed{We start with an auxiliary lemma (whose proof is due to a referee):}

\begin{Lemma}\label{L:Extra}
\changed{Let $\twocatstuff{C}$ be braided, and let 
$\twocatstuff{D}$ be symmetric. If there is a
braided functor $\functorstuff{F}\colon\twocatstuff{C}\to\twocatstuff{D}$ which is faithful, then $\twocatstuff{C}$ is symmetric.}
\end{Lemma}	

\begin{proof}
\changed{This is known, so we only give a condensed proof.
The diagram
\begin{gather*}
\begin{tikzcd}[ampersand replacement=\&]
{\functorstuff{F}(\obstuff{X})\otimes\functorstuff{F}(\obstuff{Y})} \&\& {\functorstuff{F}(\obstuff{Y})\otimes \functorstuff{F}(\obstuff{X})} \&\& {\functorstuff{F}(\obstuff{X})\otimes\functorstuff{F}(\obstuff{Y})} \\
\\
{\functorstuff{F}(\obstuff{X}\otimes\obstuff{Y})} \&\& {\functorstuff{F}(\obstuff{Y}\otimes\obstuff{X})} \&\& {\functorstuff{F}(\obstuff{X}\otimes\obstuff{Y})}
\arrow["{\morstuff{b}^{\scalebox{0.6}{$\twocatstuff{D}$}}_{\functorstuff{F}(\obstuff{X}),\functorstuff{F}(\obstuff{Y})}}", from=1-1, to=1-3]
\arrow[""{name=0, anchor=center, inner sep=0}, "{\idmor}"', from=1-1, to=1-5,bend left=25]
\arrow["{\functorstuff{F}^{2}_{\mathtt{X,Y}}}", from=1-1, to=3-1]
\arrow[""{name=1, anchor=center, inner sep=0}, "\cong"', from=1-1, to=3-1]
\arrow["{\morstuff{b}^{\scalebox{0.6}{$\twocatstuff{D}$}}_{\functorstuff{F}(\obstuff{Y}),\functorstuff{F}(\obstuff{X})}}", from=1-3, to=1-5]
\arrow["{\functorstuff{F}^{2}_{\mathtt{Y,X}}}"', from=1-3, to=3-3]
\arrow[""{name=2, anchor=center, inner sep=0}, "\cong", from=1-3, to=3-3]
\arrow["\simeq", from=1-5, to=3-5]
\arrow[""{name=3, anchor=center, inner sep=0}, "{\functorstuff{F}^{2}_{\mathtt{X,Y}}}"', from=1-5, to=3-5]
\arrow["{\functorstuff{F}(\morstuff{b}_{\mathtt{X,Y}})}"', from=3-1, to=3-3]
\arrow["{\functorstuff{F}(\morstuff{b}_{\mathtt{Y,X}})}"', from=3-3, to=3-5]
\arrow["{(2)}"{description}, draw=none, from=1, to=2]
\arrow["{(1)}"{description, pos=0.2}, draw=none, from=1-3, to=0]
\arrow["{(3)}"{description}, draw=none, from=2, to=3]
\end{tikzcd}
\end{gather*}
commutes: (1) by $\twocatstuff{D}$ being symmetric, (2) and (3) by coherence. We get that
\begin{gather*}
\End_{\scalebox{0.7}{$\twocatstuff{C}$}}(\obstuff{X}\hcirc\obstuff{Y})
\xrightarrow{\functorstuff{F}_{\obstuff{X},\obstuff{Y}}}
\End_{\scalebox{0.7}{$\twocatstuff{D}$}}\big(\functorstuff{F}(\obstuff{X}\hcirc\obstuff{Y})\big)
\xrightarrow{(\functorstuff{F}_{\obstuff{X},\obstuff{Y}}^{2})^{-1}\vcirc\placeholder\vcirc\functorstuff{F}_{\obstuff{X},\obstuff{Y}}^{2}}
\End_{\scalebox{0.7}{$\twocatstuff{D}$}}\big(\functorstuff{F}(\obstuff{X})\hcirc\functorstuff{F}(\obstuff{Y})\big)
\end{gather*}
sends both $\idmor_{\obstuff{X}\hcirc\obstuff{Y}}$ and 
$\morstuff{b}_{\obstuff{Y},\obstuff{X}}\vcirc\morstuff{b}_{\obstuff{X},\obstuff{Y}}$ to $\idmor_{\functorstuff{F}(\obstuff{X})\hcirc\functorstuff{F}(\obstuff{Y})}$. But this is a composite of
injections, and hence it is injective. Thus, we get 
$\idmor_{\obstuff{X}\hcirc\obstuff{Y}}=\morstuff{b}_{\obstuff{Y},\obstuff{X}}\vcirc\morstuff{b}_{\obstuff{X},\obstuff{Y}}$.}
\end{proof}

\changed{The flip map defines a symmetric structure on $\vect[\K]$. Hence, \autoref{L:Extra} implies the claim.}
\end{proof}

Assume that $\rank\trep<\infty$ and 
let $\addg[\trep]$ denote the additive Grothendieck group of 
the 2-representation $\trep$. (By \autoref{L:RepsSemisimple} 
we are in the semisimple case so the additive and the abelian 
Grothendieck groups agree.) For \changed{essentially fusion} categories $\twocatstuff{C}$ 
one can define $\addg[\twocatstuff{C}]$ without issue by \autoref{L:RepsKrullSchmidt} \tchanged{(even though $\twocatstuff{C}$ is allowed to have infinitely many isomorphism classes of simple objects)}.

\begin{Lemma}\label{L:RepsGrothendieck}
Write \changed{$\tworep=\trep(\obstuff{X}|\obstuff{X}\in\twocatstuff{C})$}.
The additive Grothendieck group $\addg[\tworep]$ is a 
$\addg[\twocatstuff{C}]$-representation.
\end{Lemma}

\begin{proof}
Easy to check and omitted.
\end{proof}

We write $\equirept$ for equivalence of 
$\addg[\twocatstuff{C}]$-representations.


\section{Rank one 2-representations of SL2 webs}\label{S:RepsSL2}


For the \ochanged{rest} of the paper let $\K=\C$.
As we will see, the main players in this section are complex bilinear forms.


\subsection{SL2 webs}\label{SS:RepsSL2Prelim}


We first recall the \emph{Temperley--Lieb category}, 
\ochanged{or} \emph{Rumer--Teller--Weyl category}, 
that we will call the $\sln[2]$ web category.

\begin{Definition}\label{D:RepsSL2Category}
Fix $\qpar\in\C\setminus\{0\}$.
Let $\web$ denote the $\C$-linear pivotal category $\hcirc$-generated by 
the selfdual object $\genob$, and $\vcirc$-$\hcirc$-generated by morphisms called \emph{caps} and \emph{cups} (also called \emph{bilinear form} and \emph{coform}):
\begin{gather*}
\capm=
\begin{tikzpicture}[anchorbase,scale=1]
\draw[usual] (0,0)to[out=90,in=180] (0.5,0.5)to[out=0,in=90] (1,0);
\end{tikzpicture}
\colon\genob\hcirc\genob\to\munit
,\quad
\cupm=
\begin{tikzpicture}[anchorbase,scale=1]
\draw[usual] (0,0)to[out=270,in=180] (0.5,-0.5)to[out=0,in=270] (1,0);
\end{tikzpicture}
\colon\munit\to\genob\hcirc\genob
,
\end{gather*}
modulo the $\vcirc$-$\hcirc$-ideal generated by
\emph{isotopy} and \emph{circle evaluation}:
\begin{gather*}
\begin{tikzpicture}[anchorbase,scale=1]
\draw[usual] (0,-0.5)to(0,0)to[out=90,in=180] (0.5,0.5)to[out=0,in=90] (1,0);
\draw[usual] (1,0)to[out=270,in=180] (1.5,-0.5)to[out=0,in=270] (2,0) to (2,0.5);
\end{tikzpicture}
=
\begin{tikzpicture}[anchorbase,scale=1]
\draw[usual] (0,-0.5)to(0,0.5);
\end{tikzpicture}
=
\begin{tikzpicture}[anchorbase,scale=1,xscale=-1]
\draw[usual] (0,-0.5)to(0,0)to[out=90,in=180] (0.5,0.5)to[out=0,in=90] (1,0);
\draw[usual] (1,0)to[out=270,in=180] (1.5,-0.5)to[out=0,in=270] (2,0) to (2,0.5);
\end{tikzpicture}
\,,\quad
\begin{tikzpicture}[anchorbase,scale=1]
\draw[usual] (0,0)to[out=90,in=180] (0.5,0.5)to[out=0,in=90] (1,0);
\draw[usual] (0,0)to[out=270,in=180] (0.5,-0.5)to[out=0,in=270] (1,0);
\end{tikzpicture}
=-\qnum
=-\qpar-\qpar^{-1}
.
\end{gather*}
We call $\web$ the \emph{$\sln[2]$ web category} and 
its morphism \emph{$\sln[2]$ webs}.
\end{Definition}

\begin{Remark}\label{R:RepsSL2Category}
In this \changed{and} the following sections we work over $\C$ using a ``generic'' 
$\qpar$ instead of over $\C(\qpar)$ for a variable 
$\qpar$. The situation of $\C(\qpar)$ can be discussed 
verbatim, but the linear algebra results used in this note need to be adjusted to
$\C(\qpar)$.
\end{Remark}

Let $\slnbasis{k}{l}$ denote the set of crossingless matchings
of $k$ bottom and $l$ top points, interpreted as 
$\sln[2]$ webs in the usual way.

\begin{Lemma}\label{L:RepsSL2Basis}
The set $\slnbasis{k}{l}$ is a $\C$-basis of $\Hom_{\wweb[{\sln}]}(\genob^{\hcirc k},\genob^{\hcirc l})$.
\end{Lemma}

\begin{proof}
Well-known, see e.g. \cite{Ea-tl-presentation} for a self-contained
argument that implies the claim.
\end{proof}

A \emph{nontrivial root of unity} is a 
$\qpar\notin\{1,-1\}$ with $\qpar^{k}=1$ for some $k\in\N$.

\begin{Lemma}\label{L:RepsSL2Semisimple}
We have the following.
\begin{enumerate}

\item The simple objects of $\web$ are in one-to-one correspondence 
with $\N$.

\item $\web$ is semisimple
if only if $\qpar\in\C\setminus\{0\}$ is not a nontrivial 
root of unity.

\item $\web$ is \changed{an essentially fusion} category if only if $\qpar\in\C\setminus\{0\}$ is not a nontrivial 
root of unity.

\end{enumerate}
\end{Lemma}

\begin{proof}
Recall that $\web$ can be defined integrally, meaning 
over $\Z[\qpar,\qpar^{-1}]$, and that $\web$ is integrally 
equivalent to the category of tilting modules for quantum $\sln$.
This is a type of folk theorem that dates back to \cite{RuTeWe-sl2}, 
see e.g. \cite[Theorem 2.58]{El-ladders-clasps}, \cite[Proposition 2.3]{AnStTu-semisimple-tilting} 
or \cite[Proposition 2.13]{SuTuWeZh-mixed-tilting}. 
The statements follow then from specialization to the complex numbers, 
which is well-understood on the tilting side, see e.g. 
\cite[Section 2]{AnTu-tilting}.
\end{proof}

Choose a square root $\qpar^{1/2}$ of $\qpar$. Let us define
\begin{gather}\label{Eq:RepsSL2Braiding}
\begin{tikzpicture}[anchorbase,scale=1]
\draw[usual] (1,0)to (0,1);
\draw[usual,crossline] (0,0)to (1,1);
\end{tikzpicture}
=
\qpar^{1/2}\cdot
\begin{tikzpicture}[anchorbase,scale=1]
\draw[usual] (1,0)to (1,1);
\draw[usual] (0,0)to (0,1);
\end{tikzpicture}
+\qpar^{-1/2}\cdot
\begin{tikzpicture}[anchorbase,scale=1]
\draw[usual] (0,0)to[out=90,in=180] (0.5,0.35) to[out=0,in=90] (1,0);
\draw[usual] (0,1)to[out=270,in=180] (0.5,0.65) to[out=0,in=270] (1,1);
\end{tikzpicture}
,
\begin{tikzpicture}[anchorbase,scale=1]
\draw[usual] (0,0)to (1,1);
\draw[usual,crossline] (1,0)to (0,1);
\end{tikzpicture}
=
\qpar^{-1/2}\cdot
\begin{tikzpicture}[anchorbase,scale=1]
\draw[usual] (1,0)to (1,1);
\draw[usual] (0,0)to (0,1);
\end{tikzpicture}
+\qpar^{1/2}\cdot
\begin{tikzpicture}[anchorbase,scale=1]
\draw[usual] (0,0)to[out=90,in=180] (0.5,0.35) to[out=0,in=90] (1,0);
\draw[usual] (0,1)to[out=270,in=180] (0.5,0.65) to[out=0,in=270] (1,1);
\end{tikzpicture}
\,.
\end{gather}

\begin{Lemma}\label{L:RepsSL2Braiding}
The formula \autoref{Eq:RepsSL2Braiding} endows $\web$ with the structure of 
a braided category.
\end{Lemma}

\begin{proof}
Well-known and easy to check. See also \cite[Section 2.1]{KaLi-TL-recoupling}.
\end{proof}

\begin{Notation}\label{N:RepsSL2Braiding}
As a braided category, we consider $\web$ with the structure induced by \autoref{Eq:RepsSL2Braiding}.
\end{Notation}


\subsection{The main theorem in the SL2 case}\label{SS:RepsSL2Main}


Let $\congruence$ denote \emph{matrix congruence}, that is, \tchanged{for complex $n$-by-$n$ matrices $\mat[A]$ and $\mat[B]$ we have:}
\begin{gather*}
\big(\mat[A]\congruence\mat[B]\big)\,\Leftrightarrow\;\big(\exists\mat[P]\in\matgln\colon \mat[A]=\mat[P]^{T}\mat[B]\mat[P]\big).
\end{gather*}
Note that two congruent matrices are of the same size.

\begin{Remark}\label{R:RepsSL2Congruence}
\tchanged{Recall that matrix congruence is define by ``($\mat[A]\congruence\mat[B]$) $\Leftrightarrow$ (the matrices $\mat[A]$ and $\mat[B]$ represent the same bilinear form up to change-of-basis)''.}
\end{Remark}

The proof of the following theorem is 
given in \autoref{SS:RepsSL2Proof}.

\begin{Theorem}\label{T:RepsSL2Main}
Assume $\qpar\in\C\setminus\{0\}$ is not a nontrivial 
root of unity.
\begin{enumerate}

\item Let $n\geq 2$. For every $\mat\in\matgln$ with $\trace(\mat^{T}\mat^{-1})=-\qnum$ there exists 
a simple transitive fiber 2-representation $\mrep_{\mat}^{n}$ of $\web$ constructed in the proof of \autoref{L:RepsSL2Matrix}.
\textbf{(Existence)}

\item We have $\mrep_{\mat}^{n}\equirep\mrep_{\mat[M]}^{m}$	
if and only if $\mat\congruence\mat[M]$. 
\textbf{(Non-redundant)}

\item Every simple transitive fiber 2-representation of $\web$
is of the form $\mrep_{\mat}^{n}$, and every simple transitive rank one 2-representation of $\web$ arises in this way.
\textbf{(Complete)}

\end{enumerate}
Moreover, there are infinitely many nonequivalent 
simple transitive rank one 2-representations of $\web$.
\end{Theorem}

In fact, we will make 
\autoref{T:RepsSL2Main}.(a) and (b) even more explicit.
We list some $\mrep_{\mat}^{n}$ for $n\in\{2,3\}$, while 
for $n=4$ there are infinitely many nonequivalent $\mrep_{\mat}^{n}$, 
see \autoref{L:RepsSL2List} below for details.
Moreover, \autoref{T:RepsSL2Main} and \autoref{L:RepsSL2List} together solve \autoref{C:RepsClassification} for 
$\web$.

\begin{Remark}\label{R:RepsSL2Reps}
For $n=1$ the condition $\trace(\mat^{T}\mat^{-1})=-\qnum$ becomes 
$1=-\qnum$ which has no solutions unless $\qpar\in\{\tfrac{1}{2}(-1\pm\sqrt{3})\}$.
This is the monoid case, see e.g. \cite{KhSiTu-monoidal-cryptography},
but since $\tfrac{1}{2}(-1\pm\sqrt{3})$ are nontrivial roots of unity, this 
case is not part of \autoref{T:RepsSL2Main}.
\end{Remark}

Note that \autoref{T:RepsSL2Main} shows that the classification of 
simple transitive fiber 2-representations of \ochanged{the category} $\web$ is equivalent to 
the classification of simple transitive rank one 2-representations of $\web$.
And moreover, \autoref{T:RepsSL2Main} shows that both problems 
can be considered as a subproblem of the classification of complex bilinear forms, {\cf} 
\autoref{R:RepsSL2Congruence}.
The latter has a nice known solution that we recall below.
As we will see later, see \autoref{P:ComplexitySL2}, the converse is 
also true in a precise sense. 

\changed{Here are a few bonus observations that accompany \autoref{T:RepsSL2Main}.}

\begin{Proposition}\label{P:RepsSL2Main}
\changed{We have the following.}
\begin{enumerate}

\item We have
$\addg[\mrep_{\mat}^{n}]\equirept\addg[\mrep_{\mat[M]}^{m}]$ as $\addgg[\web]$-representations if and only if $n=m$.

\item The fiber 2-representation $\mrep_{\mat}^{n}$ is braided if and only if 
($\qpar=1$, $n=2$ and $\mat$ is a standard solution as in \autoref{E:RepsSL2Standard}).

\item There exist infinitely many Hopf algebras $H$ with 
$\corep(H)\cong_{\hcirc}\web$ as monoidal categories. \changed{In particular, }infinitely many 
of these Hopf algebras are not isomorphic to $\mathcal{O}_{\qpar}\big(\matsln[2]\big)$.

\end{enumerate}
\end{Proposition}

Let us finish this section with a few (historical) remarks.

\begin{Remark}\label{R:RepsSL2HistRemark}
The category $\web$ has been around for \changed{donkey's} years 
and is a quantum version of the category constructed, albeit in 
a different language, by Rumer--Teller--Weyl \cite{RuTeWe-sl2}.
Many people have worked on this category, too many to cite here,
and it is not surprising that \autoref{T:RepsSL2Main} 
and \autoref{P:RepsSL2Main} are, in different formulations, known 
in the literature. Most prominently, \cite{Bi-repcat-qgroup} solves 
a related problem from which, after some work, one can get \autoref{T:RepsSL2Main} and \autoref{P:RepsSL2Main}.
As pointed out in \cite{Bi-repcat-qgroup}, versions of \autoref{T:RepsSL2Main} 
and \autoref{P:RepsSL2Main} are probably even older.
Having \autoref{R:RepsSame} in mind, a similar formulation also appeared in \cite{EtOs-module-categories-slq2}, see 
for example \cite[Section 3.2]{EtOs-module-categories-slq2}.
\end{Remark}

\begin{Remark}\label{R:RepsSL2HistSecond}
The case of $\web$ is one of the few web categories 
where the modular 
representation theory of the associated group is quite well-understood, 
see \cite[Section 3.4]{Do-q-schur} for a concise discussion \changed{of 
some of the main} properties. Thus, one might hope that \autoref{T:RepsSL2Main} 
generalizes to other fields than $\C$, where the 
story is not semisimple anymore. And, indeed, the paper
\cite{Os-module-cat-non-semi-slq2} has some very similar 
results. However, \autoref{R:RepsSame} does not apply 
in the nonsemisimple case.
\end{Remark}

\begin{Remark}\label{R:RepsSL2HistThird}
\autoref{P:RepsSL2Main}.(c) was used in \cite[Theorem 5.1]{CoOsTu-growth}
which in turn was the starting point of this paper.
\end{Remark}

\begin{Remark}\label{R:RepsSL2Cellular}
The category $\web$ is cellular in the sense of \cite{We-tensors-cellular-categories} or \cite{ElLa-trace-hecke}.
The same is true for the other two web categories in this paper, 
by the main result of \cite{AnStTu-cellular-tilting}
or \cite{An-tilting-cellular} and the connection to tilting modules. 
We however do not know how to use the cellular structure 
to obtain \autoref{T:RepsSL2Main} and its relatives later on.
\end{Remark}


\subsection{Proof of \autoref{T:RepsSL2Main}}\label{SS:RepsSL2Proof}


The key will be the following lemma.

\begin{Lemma}\label{L:RepsSL2Matrix}
For $n\in\Z_{\geq 2}$ let 
$\mat\in\matgln$ be a matrix satisfying $\trace(\mat^{T}\mat^{-1})=-\qnum$. 
Then there exists an associated 2-representation $\mrep$ of $\web$ on $\vect$ 
with $\dim_{\C}\mrep(\genob)=n$. Conversely, every 2-representation $\mrep$ of $\web$ on $\vect$ with $\dim_{\C}\mrep(\genob)=n$ gives such a matrix.
\end{Lemma}

\begin{proof}
Note that a 2-representation $\mrep\colon\web\to\End_{\ffincat}(\vect)$ is determined 
by specifying a $\C$-vector space $\mrep(\genob)$, a 
nondegenerate bilinear form $\mrep(\capm)$
and a nondegenerate bilinear coform $\mrep(\cupm)$ satisfying the circle evaluation \changed{and the isotopy relation.}
From a matrix $\mat$ as in the lemma we can get this data as follows.
Firstly, let $\mrep(\genob)=\C^{n}$ with fixed ordered basis 
$\{v_{1},\dots,v_{n}\}$. Writing $\mat=(m_{ij})_{1\leq i,j\leq n}$
and $\mat^{-1}=(n_{ij})_{1\leq i,j\leq n}$
in this basis we have $\trace(\mat^{T}\mat^{-1})=\sum_{1\leq i\leq n}
\sum_{1\leq j\leq n}m_{ij}n_{ij}=-\qnum$.
We then define $\mrep(\capm)$ and 
$\mrep(\cupm)$ by
\begin{gather*}
\mrep(\capm)(v_{i}\hcirc v_{j})=m_{ij},\quad
\mrep(\cupm)(1)=\sum_{1\leq i,j\leq n}n_{ij}\cdot v_{i}\hcirc v_{j}.
\end{gather*}
Since $\mat$ is invertible we get that $\mrep(\capm)$ and 
$\mrep(\cupm)$ are nondegenerate. They moreover satisfy the circle evaluation since $\trace(\mat^{T}\mat^{-1})=-\qnum$. \changed{Finally, they satisfy the isotopy relation since the coefficients $m_{ij}$ defining $\mrep(\capm)$ and the coefficients $n_{ij}$ defining $\mrep(\cupm)$ are the entries of $\mat$ and $\mat^{-1}$, respectively.}

Reading the construction backwards gives a matrix $\mat\in\matgln$ with $\trace(\mat^{T}\mat^{-1})=-\qnum$ from a 
2-representation $\mrep\colon\web\to\End_{\ffincat}(\vect)$.
\end{proof}

\begin{Example}\label{E:RepsSL2Standard}
For $x\in\C\setminus\{0\}$ we call the matrices $\mat[S]$ of the form
\begin{gather*}
\mat[S](x)=
\begin{psmallmatrix}
0 & \xpar
\\
-\qpar\xpar & 0
\end{psmallmatrix}
\text{ or }
\mat[S](x)^{\prime}=
\begin{psmallmatrix}
0 & \xpar
\\
-\qpar^{-1}\xpar & 0
\end{psmallmatrix}
\end{gather*}
the \emph{standard solutions} for $\trace(\mat^{T}\mat^{-1})=-\qnum$. 
One easily checks that $\mat[S](x)\congruence\mat[S](y)$ and 
$\mat[S](x)\congruence\mat[S](x)^{\prime}$, and we can simply 
focus on $\mat[S]=\mat[S](1)$.
\end{Example}

\begin{Lemma}\label{L:RepsSL2MatrixTwo}
For every $n\in\Z_{\geq 2}$ there exists some
$\mat\in\matgln$ with $\trace(\mat^{T}\mat^{-1})=-\qnum$.
For $n=1$ there exists no such solution.
\end{Lemma}

\begin{proof}
Let $id_{k}$ denote the $k$-by-$k$ identity matrix.
We take
\begin{gather}\label{Eq:RepsSL2MatrixExample}
\mat=
\left(\begin{array}{@{}c|cc@{}}
id_{n-2} & 0 & 0
\\
\hline
0 & 0 & 1
\\
0 & \xpar & 0
\end{array}\right)
.
\end{gather}
The matrix $\mat$ is invertible and satisfies 
$\trace(\mat^{T}\mat^{-1})=(n-2)+\xpar+\xpar^{-1}$. Thus, we can let 
$\xpar$ be a solution of $\xpar^{2}+(\qnum[2]+n-2)\xpar+1=0$ which exists since we work over $\C$.

The case $n=1$ is discussed in \autoref{R:RepsSL2Reps}.
\end{proof}

For $n\in\Z_{\geq 2}$ let us 
denote by $\mrep_{\mat}^{n}$ the 2-representation as 
constructed in the proof of \autoref{L:RepsSL2Matrix}. The existence 
is guaranteed by \autoref{L:RepsSL2MatrixTwo}. Note also that 
$\mrep_{\mat}^{n}\equirep\mrep_{\mat[M]}^{m}$ implies $n=m$ and 
\autoref{L:RepsSL2MatrixTwo} thus gives infinitely many nonequivalent rank one 2-representations
of $\web$.

\begin{Lemma}\label{L:RepsSL2MatrixEquivalence}
The 2-representation $\mrep_{\mat}^{n}$ is faithful, thus a fiber 2-representation.
\end{Lemma}

\begin{proof}
Using the basis of the hom-spaces of $\web$ given by crossingless 
matching, see \autoref{L:RepsSL2Basis}, 
\tchanged{the lemma can be proven as follows. Firstly, 
for 
$\mrep_{\mat[S]}^{2}$ this is known 
by classical results, e.g. \changed{by \cite[2. Fundamentalsatz]{RuTeWe-sl2}}, which uses the crossingless matching basis. Indeed, this references, in modern language, also shows that $\web$ is the free category generated by a nondegenerate bilinear form. This in turn, by simply copying the $n=2$ case, proves the lemma. Alternatively (and not written in 1930s German), \cite[Section 4]{Bi-repcat-qgroup} also implies the lemma.}
\end{proof}

\begin{Lemma}\label{L:RepsSL2MatrixAll}
For any fiber 2-representation
$\trep\colon\web{\to}\End_{\ffincat}(\vect)$ there exists 
a 2-representation $\mrep_{\mat}^{n}$ with $\trep\equirep\mrep_{\mat}^{n}$ 
as 2-representations of $\web$.
\end{Lemma}

\begin{proof}
From $\trep$ we can get $\mrep_{\mat}^{n}$ as follows.
View $\trep$ as a fiber functor and 
choose an ordered basis $\{v_{1},\dots,v_{n}\}$ of $\trep(\genob)$.
Then we get the lexicographically ordered basis 
$\{v_{1}\hcirc v_{1},\dots,v_{1}\hcirc v_{n},\dots,v_{n}\hcirc v_{n}\}$ of $\trep(\genob)\hcirc\trep(\genob)$. In this basis we get a $1$-by-$n^{2}$ vector $a$
determining $\trep(\capm)$ and a $n^{2}$-by-$1$ vector $b$ determining $\trep(\cupm)$.
We then rearrange $a$ 
and $b$ into $n$-by-$n$ matrices $\mat$ and $\mat^{-1}$ 
and the isotopy relation implies that these matrices, as suggested 
by their notation, are inverses. Moreover, the circle evaluation
implies that $\trace(\mat^{T}\mat^{-1})=-\qnum$.

In total, we get a 2-representation of the form $\mrep_{\mat}^{n}$. That 
$\trep\equirep\mrep_{\mat}^{n}$ holds follows by construction.
\end{proof}

From this point onward we need to assume that we are in the semisimple case.

\begin{Lemma}\label{L:RepsSL2MatrixSimple}
Assume $\qpar\in\C\setminus\{0\}$ is not a nontrivial 
root of unity. 
All rank one simple transitive 2-representations of 
$\web$ are of the form $\mrep_{\mat}^{n}$.
\end{Lemma}

\begin{proof}
By classical theory, see e.g. \cite[Chapter XII]{Tu-qgroups-3mfds} (this uses semisimplicity), we have the following property: 
Let $\catstuff{V}$ be any monoidal abelian category. Assume that 
$\obstuff{Y}\in\catstuff{V}$ has a right dual $\obstuff{Y}^{\star}$ \tchanged{ 
and there exists} an isomorphism $\morstuff{f}\colon\obstuff{Y}\to\obstuff{Y}^{\star}$ such that
\begin{gather*}
\munit\xrightarrow{ev}\obstuff{Y}\hcirc\obstuff{Y}^{\star}
\xrightarrow{\morstuff{f}\hcirc\morstuff{f}^{-1}}
\obstuff{Y}^{\star}\hcirc\obstuff{Y}\xrightarrow{coev}\munit
\end{gather*}
equals $-\qnum\cdot\morstuff{id}_{\munit}$. Then there exists a unique 
monoidal functor $\web\to\catstuff{V}$ sending $\genob$ to $\obstuff{Y}$.

Recall that \autoref{L:RepsSemisimpleRankOne} shows that for a rank one simple transitive 2-representation we can assume that $\catstuff{V}\cong\vect$, and the proof completes.
\end{proof}

\begin{Lemma}\label{L:RepsSL2Same}
Assume $\qpar\in\C\setminus\{0\}$ is not a nontrivial 
root of unity. 
Every rank one simple 
transitive 2-representations of $\web$ 
comes from a fiber 2-representation.
\end{Lemma}

\begin{proof}
We combine
\autoref{L:RepsSL2MatrixAll} and
\autoref{L:RepsSL2MatrixSimple}.
\end{proof}

\begin{Lemma}\label{L:RepsSL2MatrixEquiClasses}
Assume $\qpar\in\C\setminus\{0\}$ is not a nontrivial 
root of unity. 
We have $\mrep_{\mat}^{n}\equirep\mrep_{\mat[M]}^{m}$ 
as 2-representations of $\web$	
if and only if $\mat\congruence\mat[M]$.
\end{Lemma}

\begin{proof}
Following the same arguments as in the proof of
\autoref{L:RepsSL2MatrixSimple}, namely the characterization 
of monoidal functors $\web\to\catstuff{V}$, one obtains that 
the datum of a rank one simple transitive 2-representation 
(or, alternatively, a fiber 2-representation by \autoref{L:RepsSL2Same}) is equal 
to the datum of a $\C$-vector space and a bilinear form. 
In turn, bilinear forms are the same as matrix congruence, see 
\autoref{R:RepsSL2Congruence}, and the 
lemma follows then from the relationship of $\mat$ and its associated bilinear form.
\end{proof}

It remains to analyze matrix congruence. 
Let $\mat[J]_{n}(\lambda)$ denote an $n$-by-$n$ 
(upper triangular) Jordan block with eigenvalue 
$\lambda\in\C$. Additionally, define two matrices by
\begin{gather*}
\mat[G]_{n}=
\tchanged{\begin{psmallmatrix}
& & & & (\text{-}1)^{n}	
\\
& & & (\text{-}1)^{n\text{-}1} & (\text{-}1)^{n\text{-}1}
\\
& & \rddots & \rddots &
\\
& \text{-}1 & \text{-}1 & &
\\
1 & 1 & & &
\end{psmallmatrix}}
,\quad
\mat[H]_{2n}(\lambda)=
\left(\begin{array}{@{}c|c@{}}
0 & id_{n}
\\
\hline
\mat[J]_{n}(\lambda) & 0
\end{array}\right).
\end{gather*}
The following is a normal form 
under $\congruence$ for complex $n$-by-$n$ matrices $\mat\in\matall$:

\begin{Lemma}\label{L:RepsSL2Congruence}
Every $\mat\in\matall$ is congruent to a direct sum of matrices 
of the form $\mat[J]_{i}(0)$, $\mat[G]_{j}$ or $\mat[H]_{2k}(\lambda)$ 
with $\lambda\notin\{0,(-1)^{k+1}\}$ determined up to $\lambda\leftrightarrow\lambda^{-1}$. Moreover, for $\mat\in\matgln$ 
the matrices $\mat[J]_{i}(0)$ do not occur.
\end{Lemma}

\begin{proof}
This is \cite[Theorem 1.1]{HoSe-congruence-complex-matrices}. The tiny addition 
in the second sentence follows directly from the fact that the $\mat[J]_{i}(0)$ 
are degenerate.
\end{proof}

The matrices $\mat[G]_{n}$ and $\mat[H]_{2n}(\lambda)$ have the following associated weighted graphs with vertices labeled by the rows/columns:
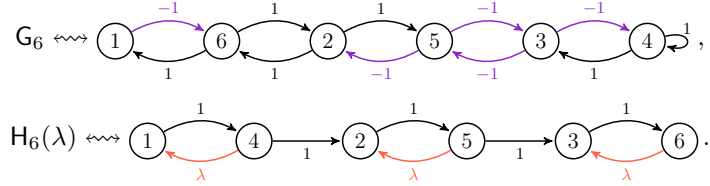
\begin{gather*}
\mat[G]_{6}\leftrightsquigarrow
\scalebox{0.7}{$\begin{tikzpicture}[anchorbase,->,>=stealth',shorten >=1pt,auto,node distance=3cm,
thick,main node/.style={circle,draw,font=\sffamily\Large\bfseries}]
\node[main node] (1) {$1$};
\node[main node] (6) [right of=1,xshift=-1.0cm] {$6$};
\node[main node] (2) [right of=6,xshift=-1.0cm] {$2$};
\node[main node] (5) [right of=2,xshift=-1.0cm] {$5$};
\node[main node] (3) [right of=5,xshift=-1.0cm] {$3$};
\node[main node] (4) [right of=3,xshift=-1.0cm] {$4$};
\path[every node/.style={font=\sffamily\small}]
(1) edge[bend left,orchid]node[above]{$-1$} (6)
(6) edge[bend left]node[above]{$1$} (2)
(2) edge[bend left]node[above]{$1$} (5)
(5) edge[bend left,orchid]node[above]{$-1$} (3)
(3) edge[bend left,orchid]node[above]{$-1$} (4)
(6) edge[bend left]node[below]{$1$} (1)
(2) edge[bend left]node[below]{$1$} (6)
(5) edge[bend left,orchid]node[below]{$-1$} (2)
(3) edge[bend left,orchid]node[below]{$-1$} (5)
(4) edge[bend left]node[below]{$1$} (3)
(4) edge [loop right]node[above]{$1$} (4);
\end{tikzpicture}$}
,\\
\mat[H]_{6}(\lambda)\leftrightsquigarrow
\scalebox{0.7}{$\begin{tikzpicture}[anchorbase,->,>=stealth',shorten >=1pt,auto,node distance=3cm,
thick,main node/.style={circle,draw,font=\sffamily\Large\bfseries}]
\node[main node] (1) {$1$};
\node[main node] (4) [right of=1,xshift=-1.0cm] {$4$};
\node[main node] (2) [right of=4,xshift=-1.0cm] {$2$};
\node[main node] (5) [right of=2,xshift=-1.0cm] {$5$};
\node[main node] (3) [right of=5,xshift=-1.0cm] {$3$};
\node[main node] (6) [right of=3,xshift=-1.0cm] {$6$};
\path[every node/.style={font=\sffamily\small}]
(1) edge[bend left]node[above]{$1$} (4)
(4) edge node[below]{$1$} (2)
(2) edge[bend left]node[above]{$1$} (5)
(5) edge node[below]{$1$} (3)
(3) edge[bend left]node[above]{$1$} (6)
(4) edge[bend left,tomato]node[below]{$\lambda$} (1)
(5) edge[bend left,tomato]node[below]{$\lambda$} (2)
(6) edge[bend left,tomato]node[below]{$\lambda$} (3);
\end{tikzpicture}$}
\,.
\end{gather*}
We display $n=4$ and $n=3$ with the general picture being similar.
Hence, the nondegenerate part of \autoref{L:RepsSL2Congruence} can be formulate 
using unions of these weighted graphs.

\begin{Example}\label{E:RepsSL2Congruence}
Let $n=2$ and take $x=1$ in \autoref{E:RepsSL2Standard}. Then 
$\mat[S]=\mat[H]_{2}(-q)$.
\end{Example}

Note that the Jordan blocks $\mat[J]_{i}(0)$ are all degenerate, so 
we can exclude them for our purposes, see the second part of \autoref{L:RepsSL2Congruence}.
For the remaining cases one directly checks that $\trace(\mat[G]_{j}^{T}\mat[G]_{j}^{-1})=(-1)^{j+1}j$
and that $\trace\big(\mat[H]_{2k}(\lambda)^{T}\mat[H]_{2k}(\lambda)^{-1}\big)=k(\lambda+\lambda^{-1}\big)$. Since $\trace(\mat^{T}\mat^{-1})$ is additive we get
\begin{gather*}
\mat[C]=\bigoplus_{a=1}^{s}\mat[G]_{j_{a}}
\oplus
\bigoplus_{b=1}^{r}\mat[H]_{2k_{b}}(\lambda_{b})
\text{ satisfies }
\trace(\mat[C]^{T}\mat[C]^{-1})
=\sum_{a=1}^{s}(-1)^{j_{a}+1}j+\sum_{b=1}^{r}k_{b}(\lambda_{b}+\lambda^{-1}_{b}).
\end{gather*}
Thus, \autoref{L:RepsSL2Congruence} gives us a list of 
solutions of $\trace(\mat^{T}\mat^{-1})=-\qnum$ up to $\congruence$. This is exactly what we want for \autoref{T:RepsSL2Main} to be as explicit as possible.

\begin{Example}\label{L:RepsSL2nTwo}
For $n=2$ we have $\trace(\mat[G]_{2}^{T}\mat[G]_{2}^{-1})=-\qnum$ 
or $\trace\big((\mat[G]_{1}\oplus\mat[G]_{1})^{T}(\mat[G]_{1}\oplus\mat[G]_{1})^{-1}\big)=-\qnum$ if and only if
$\qpar=1$ or $\qpar=-1$, while $\trace\big(\mat[H]_{2}(\lambda)^{T}\mat[H]_{2}(\lambda)^{-1}\big)=-\qnum$
if and only if $\lambda\in\{-\qpar,-\qpar^{-1}\}$.
In particular, for $\qpar\notin\{\pm 1\}$ we have $\mat[S]$ as an unique solution up to $\congruence$.
\end{Example}

\autoref{L:RepsSL2nTwo} generalizes as follows:

\begin{Lemma}\label{L:RepsSL2List}
We have the following solutions of 
$\trace(\mat^{T}\mat^{-1})=-\qnum$ up to $\congruence$.
\begin{enumerate}

\item For $n=2$ there is the solution $\mat=\mat[S]$ if $\qpar\notin\{\pm 1\}$. For $\qpar=1$ has the additional solution 
$\mat=\begin{psmallmatrix}0 & -1\\1 & 1\end{psmallmatrix}$ and 
$\qpar=-1$ has the additional solution 
$\mat=\begin{psmallmatrix}1 & 0\\0 & 1\end{psmallmatrix}$.

\item For $n=3$ there are solutions of the form
\begin{gather*}
\scalebox{0.8}{\begin{tabular}{C|C|C|C|C}
\mat & \mat[G]_{1}\oplus\mat[H]_{1}(\lambda) & \mat[G]_{1}\oplus\mat[G]_{1}\oplus\mat[G]_{1} & \mat[G]_{1}\oplus\mat[G]_{2} & \mat[G]_{3} \\
\hline
\#\text{sols} & \text{one or two} & \text{one for }\qpar\in\{\tfrac{1}{2}(-3\pm\sqrt{5})\} & \text{one for }\qpar\in\{\pm(-1)^{1/3}\} & \text{one for }\qpar\in\{\tfrac{1}{2}(-3\pm\sqrt{5})\} \\
\end{tabular}}
\end{gather*}
with $\lambda$ a root of 
$x^{2}+\tfrac{1}{2}(1+\qnum[2])x+1$ which has two solutions unless 
$\qpar\in\{\tfrac{1}{2}(3\pm\sqrt{5}),\tfrac{1}{2}(-5\pm\sqrt{21})\}$.

\item For $n=4$ there are infinitely many solutions.

\end{enumerate}
To get a complete list we use the canonical forms under orthogonal 
congruence in e.g. \cite{Ho-orthogonal-congruence}.
\end{Lemma}

\begin{proof}
Directly from the above discussion, and omitted. We only point out
two \ochanged{observations}. 

First, note that general congruence will not keep 
$\trace(\mat^{T}\mat^{-1})=-\qnum$ invariant. In particular, the 
above needs to be combined with orthogonal congruence as in the final sentence of the lemma.

Second, that 
for $n\geq 4$ we can have $\mat[H]_{k}(\lambda)\oplus\mat[H]_{l}(\mu)\oplus\text{Rest}$ 
appearing. Say $\text{Rest}$ only consists of $\mat[G]_{j}$ summands. 
Then we get infinitely 
many solutions: Fix an arbitrary $\mu$. Then the relevant equations 
for $\lambda$ always have solutions since our ground field is algebraically closed.
\end{proof}

Thus, we have proven \autoref{T:RepsSL2Main}.


\subsection{Proof of \autoref{P:RepsSL2Main}}\label{SS:RepsSL2ProofTwo}


\begin{Lemma}\label{L:RepsSL2Grothendieck}
We have $\addg[\mrep_{\mat}^{n}]\equirept\addg[\mrep_{\mat[M]}^{m}]$ 
as $\addgg[\web]$-representations \changed{if and only if $n=m$.}
\end{Lemma}

\begin{proof}
To see that we have \tchanged{$\addg[\mrep_{\mat}^{n}]\not\equirept\addg[\mrep_{\mat[M]}^{m}]$}
for $n\neq m$
we observe that 
$\addgg[\web]\cong\Z[X]$ as rings via the map 
$[\genob]\mapsto X$, and $X$ acts on 
$\addg[\mrep_{\mat}^{n}]$ by $n$.
The converse follows since the 2-representations of the form $\mrep_{\mat}^{n}$
are given by fiber functors and the twist of the bilinear form
and coform can not be detected, {\cf} \cite[Theorem 5.3.12]{EtGeNiOs-tensor-categories}.
\end{proof}

\begin{Lemma}\label{L:RepsSL2MatrixSymmetric}
The fiber 2-representation $\mrep_{\mat}^{n}$ 
is braided if and only if ($\qpar=1$, $n=2$ and $\mat$ is a standard solution as in \autoref{E:RepsSL2Standard}).
\end{Lemma}

\begin{proof}
\changed{\autoref{L:Extra} implies that $\mrep_{\mat}^{n}$ 
being braided implies that $\web$ is symmetric and that the crossing is send to the flip map.
The following calculation shows that the standard solution is the only possible choice where that happens.}

We view $\mrep_{\mat}^{n}$ as a functor $\web\to\vect$.
By the proof of \autoref{L:RepsSL2Matrix}, we have that
\begin{gather*}
\mrep_{\mat}^{n}(\cupm\vcirc\capm)(v_{i}\hcirc v_{j})
=\changed{m_{ij}\mrep_{\mat}^{n}(\cupm)(1)}=m_{ij}\sum_{1\leq k,l\leq n}n_{kl}\cdot v_{k}\hcirc v_{l}.
\end{gather*}
Hence, we get that
\begin{gather}\label{Eq:RepsSL2IsFlip}
\mrep_{\mat}^{n}
\left(
\begin{tikzpicture}[anchorbase,scale=0.5]
\draw[usual] (1,0)to (0,1);
\draw[usual,crossline] (0,0)to (1,1);
\end{tikzpicture}
\right)(v_{i}\hcirc v_{j})
=
\qpar^{1/2}\cdot v_{i}\hcirc v_{j}
+\qpar^{-1/2}\cdot
\big(
m_{ij}\sum_{1\leq k,l\leq n}n_{kl}\cdot v_{k}\hcirc v_{l}
\big)
.
\end{gather}
For this to be the flip map we then need $m_{ij}n_{ij}=-\qpar$, 
$m_{ij}n_{ji}=\qpar^{1/2}$ and $n_{kl}=0$ else. Since these have to hold for all 
$i,j\in\{1,\dots,n\}$ with $i\neq j$ we therefore need $n=2$.

For $n=2$ a direct calculation shows that the only $2$-by-$2$ matrices with 
$\trace(\mat^{T}\mat^{-1})=-\qnum$ and with \autoref{Eq:RepsSL2IsFlip} 
being the flip map are the standard solutions for $\qpar=1$.
\end{proof}

In the two latter cases in \autoref{L:RepsSL2MatrixSymmetric} the 
fiber 2-representation $\mrep_{\mat}^{n}$ is even symmetric
by \autoref{L:RepsBraidedReps}.

\begin{Example}\label{E:RepsSL2MatrixSymmetric}
We again view $\mrep_{\mat}^{n}$ as a functor $\web\to\vect$.
Let $n=3$ and take the matrix $\mat$ as in \autoref{Eq:RepsSL2MatrixExample}.
For $\qpar=1$ the variable $x$ has to be $\frac{1}{2}(-3\pm\sqrt{5})$.
For $x=\frac{1}{2}(-3+\sqrt{5})$ one gets (in an appropriate order of the basis 
$\{v_{i}\otimes v_{j}|1\leq i,j\leq n\}$) that
\begin{gather*}
\mrep_{\mat}^{3}
\left(
\begin{tikzpicture}[anchorbase,scale=1]
\draw[usual] (1,0)to (0,1);
\draw[usual,crossline] (0,0)to (1,1);
\end{tikzpicture}
\right)
=
\begin{psmallmatrix}
2 & 0 & 0 & 0 & 0 & x & 0 & 1 & 0 \\
0 & 1 & 0 & 0 & 0 & 0 & 0 & 0 & 0 \\
0 & 0 & 1 & 0 & 0 & 0 & 0 & 0 & 0 \\
0 & 0 & 0 & 1 & 0 & 0 & 0 & 0 & 0 \\
0 & 0 & 0 & 0 & 1 & 0 & 0 & 0 & 0 \\
1 & 0 & 0 & 0 & 0 & y & 0 & 1 & 0 \\
0 & 0 & 0 & 0 & 0 & 0 & 1 & 0 & 0 \\
x^{g} & 0 & 0 & 0 & 0 & 1 & 0 & y^{g} & 0 \\
0 & 0 & 0 & 0 & 0 & 0 & 0 & 0 & 1 \\
\end{psmallmatrix}
\quad\text{with}\quad
\begin{gathered}
x=\tfrac{1}{2}(-3+\sqrt{5}),
\\
x^{g}=\tfrac{1}{2}(-3-\sqrt{5}),
\\
y=\tfrac{1}{2}(-1+\sqrt{5}),
\\
y^{g}=\tfrac{1}{2}(-1-\sqrt{5}),
\end{gathered}
\end{gather*}
which squares to the identity, but is clearly not the flip map.
\end{Example}

\begin{Lemma}\label{L:RepsSL2Hopf}
\changed{There exist infinitely many Hopf algebras $H$ with 
$\corep(H)\cong\web$ as monoidal categories.
In particular, infinitely many 
of these Hopf algebras are not isomorphic to $\mathcal{O}_{\qpar}\big(\matsln[2]\big)$.}
\end{Lemma}

\begin{proof}
\changed{Reconstruction theory implies that any fiber 2-functor 
$\mrep\colon\web\to\vect$ gives rise to a Hopf algebra $H$ being 
the coend of $\mrep$. The comodules over $H$ give a monoidal category 
equivalent to $\web$.
All of this is a direct consequence 
of \cite[Theorem 4.3.1]{EtGeNiOs-tensor-categories}.
Reconstruction theory moreover implies that the resulting Hopf algebras are not isomorphic whenever the used fiber functors are not equivalent. Now we use \autoref{T:RepsSL2Main}.}
\end{proof}

The section \ochanged{is complete}.


\section{Rank one 2-representations of GL2 webs}\label{S:RepsGL2}


A lot of constructions and arguments in this section 
are similar to those in \autoref{S:RepsSL2}, so we will be brief.


\subsection{GL2 webs}\label{SS:RepsGL2Prelim}


We define webs for $\gln[2]$ as follows.

\begin{Remark}\label{R:RepsGL2Orientations}
We have two types of strands in this section with the following names:
\begin{gather*}
\text{usual}\colon
\begin{tikzpicture}[anchorbase,scale=1]
\draw[usual] (-0.5,-0.5) to (-0.5,0);
\end{tikzpicture}
\,,\quad
\text{phantom}\colon
\begin{tikzpicture}[anchorbase,scale=1]
\draw[phantom] (-0.5,-0.5) to (-0.5,0);
\end{tikzpicture}
\,.
\end{gather*}
Both types carry an orientation.
We omit the orientations in case they do not play a role in order to 
not overload the illustrations. In this case we mean any consistent orientation.
\end{Remark}

\begin{Remark}\label{R:RepsGL2Category}
Before reading \autoref{D:RepsGL2Category} 
we remind the reader that, using isotopy, one can generate
many new morphisms. For example,
\begin{gather}\label{Eq:RepsGL2Sideways}
\left(
\begin{tikzpicture}[anchorbase,scale=1]
\draw[phantom] (0.5,0)to (0.5,0.5);
\draw[usual] (0,0)to[out=90,in=180] (0.5,0.5);
\draw[usual] (1,0)to[out=90,in=0] (0.5,0.5);
\draw[usual] (1.5,0) to (1.5,0.5);
\end{tikzpicture}
\right)
\,\vcirc\,
\left(
\begin{tikzpicture}[anchorbase,scale=1]
\draw[phantom] (-0.5,-0.5) to (-0.5,0);
\draw[usual] (-1,-0.5) to (-1,0);
\draw[usual] (0,0)to[out=270,in=180] (0.25,-0.5)to[out=0,in=270] (0.5,0);
\end{tikzpicture}
\right)
=
\begin{tikzpicture}[anchorbase,scale=1]
\draw[phantom] (0.5,-0.5)to(0.5,0)to (0.5,0.5);
\draw[usual] (0,-0.5)to (0,0)to[out=90,in=180] (0.5,0.5);
\draw[usual] (1.5,0.5)to (1.5,0)to[out=270,in=0] (1.25,-0.5)to[out=180,in=270] (1,0)to[out=90,in=0] (0.5,0.5);
\end{tikzpicture}
\quad\text{is \tchanged{isotopic} to}\quad
\begin{tikzpicture}[anchorbase,scale=1]
\draw[phantom] (0.433013,-0.25)to (0,0);
\draw[usual] (0,0.5)to (0,0);
\draw[usual] (-0.433013,-0.25)to (0,0);
\end{tikzpicture}
\,.
\end{gather}
We use this silently in \autoref{D:RepsGL2Category} below.
\end{Remark}

\begin{Definition}\label{D:RepsGL2Category}
Fix $\qpar\in\C\setminus\{0\}$.
Let $\web[{\gln}]$ denote the $\C$-linear pivotal category $\hcirc$-generated by 
the \changed{dual objects} $\genob$, $\genob[Y]$, and the dual objects $\genob[P]$, $\genob[Q]$, and $\vcirc$-$\hcirc$-generated by morphisms called \emph{caps} and \emph{cups}, 
displayed and use as in \autoref{D:RepsSL2Category} but oriented:
\begin{alignat*}{2}
\capm&=
\begin{tikzpicture}[anchorbase,scale=1]
\draw[usual,directed=0.99] (0,0)to[out=90,in=180] (0.5,0.5)to[out=0,in=90] (1,0);
\end{tikzpicture}
\colon\genob\hcirc\genob[Y]\to\munit
,\quad
\cupm&=
\begin{tikzpicture}[anchorbase,scale=1]
\draw[usual,directed=0.99] (0,0)to[out=270,in=180] (0.5,-0.5)to[out=0,in=270] (1,0);
\end{tikzpicture}
\colon\munit\to\genob[Y]\hcirc\genob
,\\
\capm^{\prime}&=
\begin{tikzpicture}[anchorbase,scale=1,xscale=-1]
\draw[usual,directed=0.99] (0,0)to[out=90,in=180] (0.5,0.5)to[out=0,in=90] (1,0);
\end{tikzpicture}
\colon\genob[Y]\hcirc\genob\to\munit
,\quad
\cupm^{\prime}&=
\begin{tikzpicture}[anchorbase,scale=1,xscale=-1]
\draw[usual,directed=0.99] (0,0)to[out=270,in=180] (0.5,-0.5)to[out=0,in=270] (1,0);
\end{tikzpicture}
\colon\munit\to\genob\hcirc\genob[Y]
,
\end{alignat*}
as well 
as \emph{phantom caps} and \emph{cups}, 
\emph{phantom trilinear forms} and \emph{coforms}:
\begin{alignat*}{3}
\capp&=
\begin{tikzpicture}[anchorbase,scale=1]
\draw[phantom,directed=0.99] (0,0)to[out=90,in=180] (0.5,0.5)to[out=0,in=90] (1,0);
\end{tikzpicture}
\colon\genob[P]\hcirc\genob[Q]\to\munit
,\quad
&\cupp&=
\begin{tikzpicture}[anchorbase,scale=1]
\draw[phantom,directed=0.99] (0,0)to[out=270,in=180] (0.5,-0.5)to[out=0,in=270] (1,0);
\end{tikzpicture}
\colon\munit\to\genob[Q]\hcirc\genob[P]
,\\
\capp^{\prime}&=
\begin{tikzpicture}[anchorbase,scale=1,xscale=-1]
\draw[phantom,directed=0.99] (0,0)to[out=90,in=180] (0.5,0.5)to[out=0,in=90] (1,0);
\end{tikzpicture}
\colon\genob[Q]\hcirc\genob[P]\to\munit
,\quad
&\cupp^{\prime}&=
\begin{tikzpicture}[anchorbase,scale=1,xscale=-1]
\draw[phantom,directed=0.99] (0,0)to[out=270,in=180] (0.5,-0.5)to[out=0,in=270] (1,0);
\end{tikzpicture}
\colon\munit\to\genob[P]\hcirc\genob[Q]
,
\\
\tcapp&=
\begin{tikzpicture}[anchorbase,scale=1]
\draw[phantom,directed=0.99] (0.5,0.5)to (0.5,0);
\draw[usual,directed=0.99] (0,0)to[out=90,in=180] (0.5,0.5);
\draw[usual,directed=0.99] (1,0)to[out=90,in=0] (0.5,0.5);
\end{tikzpicture}
\colon\genob\hcirc\genob[Q]\hcirc\genob\to\munit
,\quad
&\tcupp&=
\begin{tikzpicture}[anchorbase,scale=1]
\draw[phantom,directed=0.99] (0.5,0)to (0.5,-0.5);
\draw[usual,directed=0.99] (0.5,-0.5)to[out=180,in=270] (0,0);
\draw[usual,directed=0.99] (0.5,-0.5)to[out=0,in=270] (1,0);
\end{tikzpicture}
\colon\munit\to\genob\hcirc\genob[Q]\hcirc\genob
,
\end{alignat*}
modulo the $\vcirc$-$\hcirc$-ideal generated by
\emph{isotopy} (not illustrated; we impose all possible plane isotopies), \emph{circle} and \emph{phantom circle evaluation}, \emph{H=I} and
\emph{vertical=horizontal relation} (in all consistent orientations):
\begin{gather}\label{Eq:PhantomStuff}
\begin{tikzpicture}[anchorbase,scale=1]
\draw[usual] (0,0)to[out=270,in=180] (0.5,-0.5)to[out=0,in=270] (1,0) to[out=90,in=0] (0.5,0.5)to[out=180,in=90] (0,0);
\end{tikzpicture}
=\qnum[2]
,\quad
\begin{tikzpicture}[anchorbase,scale=1]
\draw[phantom] (0,0)to[out=270,in=180] (0.5,-0.5)to[out=0,in=270] (1,0) to[out=90,in=0] (0.5,0.5)to[out=180,in=90] (0,0);
\end{tikzpicture}
=1,\quad
\begin{tikzpicture}[anchorbase,scale=1]
\draw[usual] (0,0)to (0.33,0.5)to (0,1);
\draw[phantom] (1,1)to (0.66,0.5)to (1,0);
\end{tikzpicture}
=
\begin{tikzpicture}[anchorbase,scale=1]
\draw[usual] (0,0)to (0.5,0.33);
\draw[usual] (0.5,0.66)to (0.5,0.33);
\draw[usual] (0.5,0.66)to (0,1);
\draw[phantom] (0.5,0.33)to (1,0);
\draw[phantom] (1,1) to(0.5,0.66);
\end{tikzpicture}
\,,\quad
\begin{tikzpicture}[anchorbase,scale=1]
\draw[phantom] (1,1)to (1,0);
\draw[phantom] (0,0)to (0,1);
\end{tikzpicture}
=
\begin{tikzpicture}[anchorbase,scale=1]
\draw[phantom] (0,0)to[out=90,in=180] (0.5,0.35) to[out=0,in=90] (1,0);
\draw[phantom] (1,1)to[out=270,in=0] (0.5,0.65) to[out=180,in=270] (0,1);
\end{tikzpicture}
\,.
\end{gather}
We call $\web[{\gln[2]}]$ the \emph{$\gln[2]$ web category} and 
its morphism \emph{$\gln[2]$ webs}.
\end{Definition}

\begin{Lemma}\label{L:RepsGL2Category}
We have the following.
\begin{enumerate}

\item \changed{The `oriented version' of \autoref{Eq:RepsSL2Braiding} given by e.g.
\begin{gather*}
\begin{tikzpicture}[anchorbase,scale=1]
\draw[usual,directed=0.99] (1,0)to (0,1);
\draw[usual,crossline,directed=0.99] (0,0)to (1,1);
\end{tikzpicture}
=
\qpar^{1/2}\cdot
\begin{tikzpicture}[anchorbase,scale=1]
\draw[usual,directed=0.99] (1,0)to (1,1);
\draw[usual,directed=0.99] (0,0)to (0,1);
\end{tikzpicture}
+\qpar^{-1/2}\cdot
\begin{tikzpicture}[anchorbase,scale=1]
\draw[phantom,directed=0.99] (0.5,0.3)to (0.5,0.7);
\draw[usual] (0,0)to[out=90,in=180] (0.5,0.3);
\draw[usual,directed=0.99] (0.5,0.7) to[out=180,in=270] (0,0.98) to (0,1);
\draw[usual] (0.5,0.3) to[out=0,in=90] (1,0);
\draw[usual,directed=0.99] (0.5,0.7) to[out=0,in=270] (1,0.98) to (1,1);
\end{tikzpicture}
\end{gather*}
and additionally}
\begin{gather*}
\begin{tikzpicture}[anchorbase,scale=1]
\draw[phantom,directed=0.99] (1,0)to (0,1);
\draw[usual] (0,0)to (1,1);
\end{tikzpicture}
=
\begin{tikzpicture}[anchorbase,scale=1]
\draw[phantom,directed=0.99] (1,0)to (0.75,0.75);
\draw[phantom,directed=0.99] (0.25,0.25)to (0,1);
\draw[usual] (0,0)to (1,1);
\end{tikzpicture}
\,,\quad
\begin{tikzpicture}[anchorbase,scale=1]
\draw[phantom,directed=0.99] (1,0)to (0,1);
\draw[phantom,directed=0.99] (0,0)to (1,1);
\end{tikzpicture}
=
\begin{tikzpicture}[anchorbase,scale=1]
\draw[phantom,directed=0.99] (1,0)to (1,1);
\draw[phantom,directed=0.99] (0,0)to (0,1);
\end{tikzpicture}
\,,
\end{gather*}
and similar formulas define a braiding on $\web[{\gln[2]}]$ with the 
phantom strings being part of a symmetric structure where the Reidemeister I 
relations holds \tchanged{(the full subcategory generated by $\genob[P]$ and $\genob[Q]$ is symmetric with the phantom crossing)}.

\item We have the \emph{trilinear evaluation}:
\begin{gather*}
\begin{tikzpicture}[anchorbase,scale=1]
\draw[phantom] (0.5,-0.5) to (0.5,0.5);
\draw[usual] (0,0)to[out=270,in=180] (0.5,-0.5)to[out=0,in=270] (1,0) to[out=90,in=0] (0.5,0.5)to[out=180,in=90] (0,0);
\end{tikzpicture}
=\qnum[2].
\end{gather*}

\end{enumerate}
\end{Lemma}

\begin{proof}
\textit{(a)---non-mixed part.} Easy and we just give one calculation:
\begin{gather*}
\begin{tikzpicture}[anchorbase,scale=1,rounded corners]
\draw[phantom,directed=0.99] (1,1)to (0,0.5)to (1,0);
\draw[phantom,directed=0.99] (0,0)to (1,0.5)to (0,1);
\end{tikzpicture}
=
\begin{tikzpicture}[anchorbase,scale=1]
\draw[phantom,directed=0.99] (0,0)to[out=90,in=180] (0.5,0.35) to[out=0,in=90] (1,0);
\draw[phantom,directed=0.99] (1,1)to[out=270,in=0] (0.5,0.65) to[out=180,in=270] (0,1);
\draw[phantom,directed=0.99] (-0.75,0.5)to[out=270,in=180] (-0.5,0.25)to[out=0,in=270] (-0.25,0.5) to[out=90,in=0] (-0.5,0.75)to[out=180,in=90] (-0.75,0.5);
\end{tikzpicture}
=
\begin{tikzpicture}[anchorbase,scale=1]
\draw[phantom,directed=0.99] (0,0)to[out=90,in=180] (0.5,0.35) to[out=0,in=90] (1,0);
\draw[phantom,directed=0.99] (1,1)to[out=270,in=0] (0.5,0.65) to[out=180,in=270] (0,1);
\end{tikzpicture}
=
\begin{tikzpicture}[anchorbase,scale=1,rounded corners]
\draw[phantom,directed=0.99] (1,1)to (1,0);
\draw[phantom,directed=0.99] (0,0)to (0,1);
\end{tikzpicture}
\,.
\end{gather*}
This uses the phantom circle evaluation and vertical=horizontal relation.

\textit{(a)---mixed part.} One first shows that
\begin{gather}\label{Eq:RepsGL2Digon}
\begin{tikzpicture}[anchorbase,scale=1]
\draw[phantom,directed=0.99] (0,0.2) to[out=0,in=270] (0.5,0.5) to[out=90,in=0] (0,0.8);
\draw[usual] (0,0) to (0,1);
\end{tikzpicture}
=
\begin{tikzpicture}[anchorbase,scale=1]
\draw[phantom,directed=0.99] (0,0.8) to[out=0,in=90] (0.5,0.5) to[out=270,in=0] (0,0.2);
\draw[usual] (0,0) to (0,1);
\end{tikzpicture}
=
\begin{tikzpicture}[anchorbase,scale=1]
\draw[usual] (0,0) to (0,1);
\end{tikzpicture}
=
\begin{tikzpicture}[anchorbase,scale=1,xscale=-1]
\draw[phantom,directed=0.99] (0,0.8) to[out=0,in=90] (0.5,0.5) to[out=270,in=0] (0,0.2);
\draw[usual] (0,0) to (0,1);
\end{tikzpicture}
=
\begin{tikzpicture}[anchorbase,scale=1,xscale=-1]
\draw[phantom,directed=0.99] (0,0.2) to[out=0,in=270] (0.5,0.5) to[out=90,in=0] (0,0.8);
\draw[usual] (0,0) to (0,1);
\end{tikzpicture}
\,.
\end{gather}
This is a direct consequence of the vertical=horizontal relation. Using this and similar formulas, one can show 
that the above defines a braiding.

\changed{\textit{(b).}} Immediately from \autoref{Eq:RepsGL2Digon}.
\end{proof}

For two objects $\obstuff{A},\obstuff{B}
\in\web[{\gln}]$ let $\glnbasis{A}{B}$ denote any 
(fixed) choice of 
placement of phantom edges such that the $\gln$ web obtained 
by removing the phantom edges corresponds to a crossingless matching.

\begin{Lemma}\label{L:RepsGL2Basis}
The set $\glnbasis{A}{B}$ is a $\C$-basis of $\Hom_{\wweb[{\gln}]}(\obstuff{A},\obstuff{B})$.
\end{Lemma}

\begin{proof}
Directly by using the braiding in \autoref{L:RepsGL2Category} and the 
usual crossingless matching basis of $\web$, see \autoref{L:RepsSL2Basis}. \tchanged{In more details, the relations involving phantom strings ensure that we have two cases. Firstly, a phantom string touches a usual string an even number of times. Then the phantom string can be unplugged from the usual string. On the other hand, if they touch an odd number of times, then the phantom string can be unplugged up to one attachment, and this attachment can be placed arbitrarily along the usual string. This in turn implies that the usual crossingless matching basis plus an arbitrary, but fixed and minimal, placement of phantom strings gives a basis.}
\end{proof}

The remainder of \autoref{SS:RepsSL2Prelim} goes through for $\web[{\gln[2]}]$ 
with one mild change\changed{, namely \autoref{L:RepsGL2Semisimple}.(a)}. That is:

\begin{Lemma}\label{L:RepsGL2Semisimple}
We have the following.
\begin{enumerate}

\item The simple objects of $\web[{\gln}]$ are in one-to-one correspondence 
with $\N\times\Z$.

\item $\web[{\gln}]$ is semisimple
if only if $\qpar\in\C\setminus\{0\}$ is not a nontrivial 
root of unity.

\item $\web[{\gln}]$ is \changed{an essentially fusion} category if only if $\qpar\in\C\setminus\{0\}$ is not a nontrivial 
root of unity.

\end{enumerate}
\end{Lemma}

\begin{proof}
The statement follows similarly as in \autoref{L:RepsSL2Semisimple}.
\end{proof}


\subsection{The main theorem in the GL2 case}\label{SS:RepsGL2Main}


For us an \emph{third order tensor} is an $l$-by-$m$-by-$n$ array of complex numbers. We represent a third order tensor by \changed{$\ten=(t_{ijk})_{i,j,k}$}
with $t_{ijk}\in\C$. The indexes are the rows and columns, 
as for usual matrices, and the \emph{pages} $k$.
Here is an illustration for $l=m=n=3$:
\begin{gather}\label{Eq:RepsSO3Trilinear}
\begin{tikzpicture}[anchorbase,every node/.style={anchor=north east,fill=white,minimum width=1.0cm,minimum height=5mm}]
\matrix (mA) [draw,matrix of math nodes,ampersand replacement=\&]
{
t_{{\color{spinach}1}{\color{blue}1}{\color{orchid}3}} \& t_{{\color{spinach}1}{\color{blue}2}{\color{orchid}3}} \& t_{{\color{spinach}1}{\color{blue}3}{\color{orchid}3}} \\
t_{{\color{spinach}2}{\color{blue}1}{\color{orchid}3}} \& t_{{\color{spinach}2}{\color{blue}2}{\color{orchid}3}} \& t_{{\color{spinach}2}{\color{blue}3}{\color{orchid}3}} \\
t_{{\color{spinach}3}{\color{blue}1}{\color{orchid}3}} \& t_{{\color{spinach}3}{\color{blue}2}{\color{orchid}3}} \& t_{{\color{spinach}3}{\color{blue}3}{\color{orchid}3}} \\
};
\matrix (mB) [draw,matrix of math nodes,ampersand replacement=\&] at ($(mA.south west)+(1.5,0.7)$)
{
t_{{\color{spinach}1}{\color{blue}1}{\color{orchid}2}} \& t_{{\color{spinach}1}{\color{blue}2}{\color{orchid}2}} \& t_{{\color{spinach}1}{\color{blue}3}{\color{orchid}2}} \\
t_{{\color{spinach}2}{\color{blue}1}{\color{orchid}2}} \& t_{{\color{spinach}2}{\color{blue}2}{\color{orchid}2}} \& t_{{\color{spinach}2}{\color{blue}3}{\color{orchid}2}} \\
t_{{\color{spinach}3}{\color{blue}1}{\color{orchid}2}} \& t_{{\color{spinach}3}{\color{blue}2}{\color{orchid}2}} \& t_{{\color{spinach}3}{\color{blue}3}{\color{orchid}2}} \\
};
\matrix (mC) [draw,matrix of math nodes,ampersand replacement=\&] at ($(mB.south west)+(1.5,0.7)$)
{
t_{{\color{spinach}1}{\color{blue}1}{\color{orchid}1}} \& t_{{\color{spinach}1}{\color{blue}2}{\color{orchid}1}} \& t_{{\color{spinach}1}{\color{blue}3}{\color{orchid}1}} \\
t_{{\color{spinach}2}{\color{blue}1}{\color{orchid}1}} \& t_{{\color{spinach}2}{\color{blue}2}{\color{orchid}1}} \& t_{{\color{spinach}2}{\color{blue}3}{\color{orchid}1}} \\
t_{{\color{spinach}3}{\color{blue}1}{\color{orchid}1}} \& t_{{\color{spinach}3}{\color{blue}2}{\color{orchid}1}} \& t_{{\color{spinach}3}{\color{blue}3}{\color{orchid}1}} \\
};
\draw[dashed](mA.north east)--(mC.north east);
\draw[dashed](mA.north west)-- node[sloped,above] {\small {\color{orchid}Page}} (mC.north west);
\draw[dashed](mA.south east)--(mC.south east);
\node [above left,opacity=0,text opacity=1] at ($(mC.north)+(0.25,0)$) {\small {\color{blue}Column}};
\node [left] at (mC.west) {\small {\color{spinach}Row}};
\end{tikzpicture}
\;.
\end{gather}
As before, fix bases $\{v_{1},\dots,v_{a}\}$ of $\C^{a}$.
It is immediate that a third order tensor gives a trilinear form by
\begin{gather*}
\ten\colon
\C^{l}\hcirc\C^{m}\hcirc\C^{n}\to\C,
v_{i}\hcirc v_{j}\hcirc v_{k}\mapsto t_{ijk}
.
\end{gather*}
In other words, $\ten$ is a $1$-by-$lmn$ matrix.

Assume that we have already fixed a matrix $\mat$ that corresponds 
to a bilinear form. We therefore get matrices associated to caps and cups. Let us call these 
matrices $\mat(cap)$ and $\mat(cup)$, respectively. \changed{Define
\begin{gather*}
\ten_{l}=(\ten\hcirc id_{\C^{n}})\vcirc\big(id_{\C^{l}}\hcirc id_{\C^{m}}\hcirc\mat(cup)\big)
,\\
\ten^{l}=\big(id_{\C^{l}}\hcirc id_{\C^{m}}\hcirc\mat(cap)\big)\vcirc
(\ten^{\prime}\hcirc id_{\C^{n}}),
\end{gather*}
where $\ten^{\prime}$ is the transpose tensor.}
The picture to keep in mind is \autoref{Eq:RepsGL2Sideways} 
which displays the diagrammatic 
interpretation of $\ten_{l}$.

Let us denote the set of $l$-by-$m$-by-$n$ tensors by $\tensors{l}{m}{n}$, and 
for elements in that set let us write $\ten\congruence\ten[U]$ for \emph{congruence of 
third order tensors} \tchanged{in the sense of e.g. \cite[Section 4]{BeSe-matrix-problems}, meaning, roughly speaking,} that they define the same trilinear form up to change-of-basis.
\changed{Below we write $\mat[P](cap)=\mat[P](cup)=\mat[P]$ to highlight how the next display fits to \autoref{Eq:PhantomStuff}.}

\begin{Theorem}\label{T:RepsGL2Main}
Assume $\qpar\in\C\setminus\{0\}$ is not a nontrivial 
root of unity.
\begin{enumerate}

\item Let $n\geq 2$. For every triple $\bN=(\mat,\mat[P],\ten)\in
\matgln\times\{\pm 1\}\times\tensors{n}{1}{n}$
\changed{with 
\begin{gather*}
\trace(\mat^{T}\mat^{-1})=\qnum,\quad
id_{\C^{l}\hcirc\C}=\ten^{l}\ten_{l},\quad
id_{\C\hcirc\C}=\mat[P](cap)\mat[P](cup)
\end{gather*}
there} exists 
a simple transitive fiber 2-representation $\mrep_{\bN}$ of $\web[{\gln}]$ constructed in the proof of \autoref{L:RepsGL2Matrix}.
\textbf{(Existence)}

\item We have $\mrep_{(\mat,\mat[P],\ten)}\equirep\mrep_{(\mat[M],\mat[Q],\ten[U])}$	
if and only if $\mat\congruence\mat[M]$, $\mat[P]=\mat[Q]$ 
and $\ten\congruence\mat[U]$. 
\textbf{(Non-redundant)}

\item Every simple transitive fiber 2-representation of $\web[{\gln}]$
is of the form $\mrep_{\mat}^{n}$, and every simple transitive rank one 2-representation of $\web[{\gln}]$ arises in this way.
\textbf{(Complete)}

\end{enumerate}
Moreover, there are infinitely many nonequivalent 
simple transitive rank one 2-representations of $\web[{\gln}]$.
\end{Theorem}

As before, we list some $\mrep_{\bN}$ for $n\in\{2,3\}$, while 
for $n=4$ there are infinitely many nonequivalent $\mrep_{\bN}$, 
see \autoref{L:RepsGL2List} below for details (also as before, \autoref{T:RepsGL2Main} and \autoref{L:RepsGL2List} taken together solve \autoref{C:RepsClassification} for $\web[{\gln}]$),
and the proof of \autoref{T:RepsGL2Main} gets its own section.

Moreover, we leave it to the reader to spell out the 
$\gln$ analog of \autoref{P:RepsSL2Main} (which reads essentially the same).
We rather wrap-up this section with a (historical) remark and another remark:

\begin{Remark}\label{R:RepsGL2HistRemark}
$\web[{\gln}]$ was first considered to construct a functorial version of 
Khovanov homology \cite{Bl-gl2-foams}, and $\gln$ webs have been \ochanged{studied 
intensively} in the context of link homologies, see 
e.g. \cite{EhStTu-blanchet-khovanov}, \cite{EhStTu-gl2foams}, \cite{BeHoPuWe-sl2gl2} or \cite{KrWe-gl2-homotopy-type}.
\ochanged{Indeed}, our presentation of $\web[{\gln}]$ is stolen 
from \cite{BeHoPuWe-sl2gl2}.
\end{Remark}

\begin{Remark}\label{R:RepsGL2HistRemarkTwo}
The reader familiar with \cite{Mr-quantum-gl2} will notice that 
the main theorem of that paper and \autoref{T:RepsGL2Main} are different.
This is due to us using diagrammatics 
that are not used in \cite{Mr-quantum-gl2}.
Hence, \autoref{T:RepsGL2Main} appears to be new in the presented form, 
and is the expected extension of \autoref{T:RepsSL2Main}.
\end{Remark}


\subsection{Proof of \autoref{T:RepsGL2Main}}\label{SS:RepsGL2Proof}


The proof 
of \autoref{T:RepsGL2Main} is, of course, similar to the proof 
of \autoref{T:RepsSL2Main} so we will be rather brief and focus on
the main differences.

Let $\peb\subset\web[{\gln}]$ denote the full subcategory 
generated by $\genob[P],\genob[Q]$.

\begin{Lemma}\label{L:RepsGL2DimOne}
We have the following.
\begin{enumerate}

\item For $\mat[P]\in\{\pm 1\}$ 
there exists 
a simple transitive fiber 2-representation $\mrep_{\mat[P]}$ of $\peb$ constructed similarly as in the proof of \autoref{L:RepsSL2Matrix}.
\textbf{(Existence)}

\item $\mrep_{+1}$ is not equivalent to $\mrep_{-1}$ as 2-representations of $\peb$.
\textbf{(Non-redundant)}

\item Every simple transitive fiber 2-representation of $\peb$
is of the form $\mrep_{\pm 1}$, and every simple transitive rank one 2-representation of $\peb$ arises in this way.
\textbf{(Complete)}

\end{enumerate}
\end{Lemma}

\begin{proof}
The proof is similar, but much easier, than the proof of 
\autoref{T:RepsSL2Main}. So let us only give the two new observations
needed for the proof.

Assume that we have a one column $cap=(a_{1},\dots,a_{n})^{T}$ and a one row matrix $cup=(b_{1},\dots,b_{n})$. Then
\begin{gather*}
cap\times cup=(a_{1}b_{1}+\dots+a_{n}b_{n})
,\quad	
\text{the diagonal of }cup\times cap\text{ is }(a_{1}b_{1},\dots,a_{n}b_{n}).
\end{gather*}
In particular, $cap\times cup=(1)$ and $cup\times cap=id_{n}$ can only 
hold for $n=1$. Moreover, for $n=1$ the only possible solutions are 
$a_{1}=b_{1}=\pm 1$. Thus, a 2-representation of $\peb$ 
needs to send both generating objects to $\C$, and the phantom caps
and cups to multiplication by $\pm 1$.

It then follows from the phantom circle removal and the isotopy relations 
that fixing $\pm 1$ as the value for $\capp$ determines the other three 
bilinear (co)forms, so we only have $\pm 1$ to vary.
That is:
\begin{gather*}
\left(
\begin{tikzpicture}[anchorbase,scale=1]
\draw[phantom,directed=0.99] (0,0)to[out=90,in=180] (0.5,0.5)to[out=0,in=90] (1,0);
\end{tikzpicture}
\mapsto\cdot-1\right)
\Rightarrow
\left(
\begin{tikzpicture}[anchorbase,scale=1,xscale=-1]
\draw[phantom,directed=0.99] (0,0)to[out=270,in=180] (0.5,-0.5)to[out=0,in=270] (1,0);
\end{tikzpicture}
\mapsto\cdot-1\right)
\quad\text{since}\quad
\begin{tikzpicture}[anchorbase,scale=1,xscale=-1]
\draw[phantom,directed=0.99] (0,0)to[out=270,in=180] (0.5,-0.5)to[out=0,in=270] (1,0) to[out=90,in=0] (0.5,0.5)to[out=180,in=90] (0,0);
\end{tikzpicture}
=1,
\\
\left(
\begin{tikzpicture}[anchorbase,scale=1]
\draw[phantom,directed=0.99] (0,0)to[out=90,in=180] (0.5,0.5)to[out=0,in=90] (1,0);
\end{tikzpicture}
\mapsto\cdot-1\right)
\Rightarrow
\left(
\begin{tikzpicture}[anchorbase,scale=1,xscale=1]
\draw[phantom,directed=0.99] (0,0)to[out=270,in=180] (0.5,-0.5)to[out=0,in=270] (1,0);
\end{tikzpicture}
\mapsto\cdot-1\right)
\quad\text{since}\quad
\begin{tikzpicture}[anchorbase,scale=1]
\draw[phantom,directed=0.99] (0,-0.5)to (0,0)to[out=90,in=180] (0.5,0.5)to[out=0,in=90] (1,0)to[out=270,in=180] (1.5,-0.5)to[out=0,in=270] (2,0) to (2,0.5);	
\end{tikzpicture}
=
\begin{tikzpicture}[anchorbase,scale=1]
\draw[phantom,directed=0.99] (0,-0.5)to (0,0.5);	
\end{tikzpicture}
\,,
\\
\left(
\begin{tikzpicture}[anchorbase,scale=1,xscale=-1]
\draw[phantom,directed=0.99] (0,0)to[out=270,in=180] (0.5,-0.5)to[out=0,in=270] (1,0);
\end{tikzpicture}
\mapsto\cdot-1\right)
\Rightarrow
\left(
\begin{tikzpicture}[anchorbase,scale=1,xscale=-1]
\draw[phantom,directed=0.99] (0,0)to[out=90,in=180] (0.5,0.5)to[out=0,in=90] (1,0);
\end{tikzpicture}
\mapsto\cdot-1\right)
\quad\text{since}\quad
\begin{tikzpicture}[anchorbase,scale=1,xscale=-1]
\draw[phantom,directed=0.99] (0,-0.5)to (0,0)to[out=90,in=180] (0.5,0.5)to[out=0,in=90] (1,0)to[out=270,in=180] (1.5,-0.5)to[out=0,in=270] (2,0) to (2,0.5);	
\end{tikzpicture}
=
\begin{tikzpicture}[anchorbase,scale=1,xscale=-1]
\draw[phantom,directed=0.99] (0,-0.5)to (0,0.5);	
\end{tikzpicture}
\,
.
\end{gather*}
All other cases follow by symmetry.
\end{proof}

Similarly as above, let $\qeb\subset\web[{\gln}]$ denote the full subcategory 
generated by $\genob,\genob[Y]$.

\begin{Lemma}\label{L:RepsGL2DimOneTwo}
\autoref{T:RepsSL2Main} holds verbatim for $\qeb$.
\end{Lemma}

\begin{proof}
As in the proof of \autoref{L:RepsGL2DimOne},
\begin{gather*}
\left(
\begin{tikzpicture}[anchorbase,scale=1]
\draw[usual,directed=0.99] (0,0)to[out=90,in=180] (0.5,0.5)to[out=0,in=90] (1,0);
\end{tikzpicture}
\mapsto\text{fixed}\right)
\Rightarrow
\left(
\begin{tikzpicture}[anchorbase,scale=1,xscale=1]
\draw[usual,directed=0.99] (0,0)to[out=270,in=180] (0.5,-0.5)to[out=0,in=270] (1,0);
\end{tikzpicture}
\mapsto\text{fixed}\right)
\\
\text{since}\quad
\begin{tikzpicture}[anchorbase,scale=1]
\draw[usual,directed=0.99] (0,-0.5)to (0,0)to[out=90,in=180] (0.5,0.5)to[out=0,in=90] (1,0)to[out=270,in=180] (1.5,-0.5)to[out=0,in=270] (2,0) to (2,0.5);	
\end{tikzpicture}
=
\begin{tikzpicture}[anchorbase,scale=1]
\draw[usual,directed=0.99] (0,-0.5)to (0,0.5);	
\end{tikzpicture}
\quad\text{and}\quad
\begin{tikzpicture}[anchorbase,scale=1,yscale=-1]
\draw[usual,directed=0.99] (0,-0.5)to (0,0)to[out=90,in=180] (0.5,0.5)to[out=0,in=90] (1,0)to[out=270,in=180] (1.5,-0.5)to[out=0,in=270] (2,0) to (2,0.5);	
\end{tikzpicture}
=
\begin{tikzpicture}[anchorbase,scale=1,yscale=-1]
\draw[usual,directed=0.99] (0,-0.5)to (0,0.5);	
\end{tikzpicture}
\,,
\end{gather*}
etc. \tchanged{(as above, the cup oriented rightwards and the circle evaluation fixes the assignment for the cap oriented leftwards, and then the zigzag fixes the assignment for the cup oriented leftwards).} The rest of the proof works, mutatis mutandis, as for $\sln$.
\end{proof}

\changed{We will refer to the triples $\bN=(\mat,\mat[P],\ten)$ in \autoref{T:RepsGL2Main} as \emph{$\gln$ triples}.}

\begin{Lemma}\label{L:RepsGL2Matrix}
For $n\in\Z_{\geq 2}$ let 
$\bN=(\mat,\mat[P],\ten)$ be a $\gln$ triple.
Then there exists an associated 2-representation $\mrep$ of $\web[{\gln}]$ on $\vect$ 
with $\dim_{\C}\mrep(\genob)=\dim_{\C}\mrep(\genob[Y])=n$ and 
$\dim_{\C}\mrep(\genob[P])=\dim_{\C}\mrep(\genob[Q])=1$. Conversely, every 2-representation $\mrep$ of $\web[{\gln}]$ on $\vect$ with $\dim_{\C}\mrep(\genob[P])=\dim_{\C}\mrep(\genob[Q])=1$ gives such a triple.
\end{Lemma}

\begin{proof}
Very similar to the proof of \autoref{L:RepsGL2Matrix} 
with the following two differences. Firstly, the phantom part is taken care of 
by \autoref{L:RepsGL2DimOne} while the 
usual part is \autoref{L:RepsGL2DimOneTwo}. 
The two sides are related via the trilinear form and the H=I relation.
\changed{Note that \autoref{L:RepsGL2Category} shows that the trilinear form determines 
the trilinear coform in exactly the same way as the bilinear 
form and coform are related, so we only need to specify the trilinear form. }
Finally, the H=I relation is part of the definition.
\end{proof}

\begin{Lemma}\label{L:RepsGL2MatrixTwo}
For every $n\in\Z_{\geq 2}$ there exists some $\gln$ triple.
For $n=1$ there exists no such triple.
\end{Lemma}

\begin{proof}
Let us take $\mat[P]=1$, and let $\mat$ be any matrix satisfying 
$\trace(\mat^{T}\mat^{-1})=\qnum$. The existence of the latter 
is guaranteed by (the same arguments as in) \autoref{L:RepsSL2MatrixTwo}, while the choice $\mat[P]=1$ satisfies $\trace(\mat[P]^{T}\mat[P]^{-1})=1$
and $id=\mat[P](cap)\mat[P](cup)$.
\changed{We may construct a trilinear form $\ten$ by mapping $v_{i}\hcirc 1\hcirc v_{k}\mapsto\mat(cap)(v_{i}\hcirc v_{k})$. This form satisfies the conditions of
\autoref{T:RepsGL2Main}. On the diagrammatic side, this corresponds to ignoring the
phantom edges and identifying the bilinear form with the trilinear form.
The property 
$id=\ten^{l}\ten_{l}$ is then clearly satisfied.}

For $n=1$ see \autoref{R:RepsSL2Reps}.
\end{proof}

As before, let $\mrep_{\bN}$ be the 2-representation constructed 
above. 

\begin{Lemma}\label{L:RepsGL2MatrixEquivalence}
The 2-representation $\mrep_{\bN}$ is faithful, thus a fiber 2-representation.
\end{Lemma}

\begin{proof}
As before, but using \autoref{L:RepsGL2Basis} instead of \autoref{L:RepsSL2Basis}.
\end{proof}

\begin{Lemma}\label{L:RepsGL2Collect}
The statements \autoref{L:RepsSL2MatrixAll} to 
\autoref{L:RepsSL2MatrixEquiClasses}
hold mutatis mutandis for $\web[{\gln}]$ as well.
\end{Lemma}

\begin{proof}
Only two things changes with respect to the proofs given in 
\autoref{SS:RepsSL2Proof}. Firstly, one uses the basis in
\autoref{L:RepsGL2Basis} instead of the crossingless matching basis.
Second, the careful 
reader can copy the arguments in \cite[Chapter XII]{Tu-qgroups-3mfds} 
to get the analog of the result used in the proof of \autoref{L:RepsSL2MatrixSimple}.
\end{proof}

\begin{Lemma}\label{L:RepsGL2List}
$\gln$ triples, up to $\congruence$, are given by:
\begin{enumerate}[label=$\triangleright$]

\item The matrix $\mat$ is classified as in \autoref{L:RepsSL2List}.

\item The sign $\mat[P]$ can be chosen freely.

\item The tensor $\ten$ is classified as $\mat$ in \autoref{L:RepsSL2List} together with the choice of a sign.

\end{enumerate}
\end{Lemma}

\begin{proof}
The only extra information one needs beyond \autoref{L:RepsSL2List} 
is the classification of $n$-by-$1$-by-$n$ trilinear forms, which is 
the same as the classification of $n$-by-$n$ bilinear forms 
up to a sign. This is easy to see, but can also be found 
explicitly spelled out in \cite[Introduction]{Th-trilinear}.
\end{proof}

Hence, taking the above together proves \autoref{T:RepsGL2Main}.


\section{Rank one 2-representations of SO3 webs}\label{S:RepsSO3}


As expected, a lot of constructions and arguments in this section are 
similar to those in the previous sections, so we will be brief 
and focus on the new bits.


\subsection{SO3 webs}\label{SS:RepsSO3Prelim}


We start with a reminder on the $\son[3]$ web category.
As in the previous section we
silently use (an analog of) \autoref{R:RepsGL2Category}.

\begin{Definition}\label{D:RepsSO3Category}
Fix $\qpar\in\C\setminus\{0\}$ with $\qpar^{2}+\qpar^{-2}\neq 0$.
Let $\web[{\son[3]}]$ denote the $\C$-linear pivotal category $\hcirc$-generated by 
the selfdual object $\genob$, and $\vcirc$-$\hcirc$-generated by morphisms called \emph{bilinear} and \emph{trilinear forms} and \emph{coforms}:
\begin{alignat*}{3}
\capm&=
\begin{tikzpicture}[anchorbase,scale=1]
\draw[usual] (0,0)to[out=90,in=180] (0.5,0.5)to[out=0,in=90] (1,0);
\end{tikzpicture}
\colon\genob\hcirc\genob\to\munit
,\quad
&\cupm&=
\begin{tikzpicture}[anchorbase,scale=1]
\draw[usual] (0,0)to[out=270,in=180] (0.5,-0.5)to[out=0,in=270] (1,0);
\end{tikzpicture}
\colon\munit\to\genob\hcirc\genob
,\\
\tcapm&=
\begin{tikzpicture}[anchorbase,scale=1]
\draw[usual] (0,0)to[out=90,in=180] (0.5,0.5);
\draw[usual] (1,0)to[out=90,in=0] (0.5,0.5);
\draw[usual] (0.5,0)to (0.5,0.5);
\end{tikzpicture}
\colon\genob\hcirc\genob\hcirc\genob\to\munit
,\quad
&\tcupm&=
\begin{tikzpicture}[anchorbase,scale=1]
\draw[usual] (0.5,-0.5)to[out=180,in=270] (0,0);
\draw[usual] (0.5,-0.5)to[out=0,in=270] (1,0);
\draw[usual] (0.5,-0.5)to (0.5,0);
\end{tikzpicture}
\colon\munit\to\genob\hcirc\genob\hcirc\genob
,
\end{alignat*}
modulo the $\vcirc$-$\hcirc$-ideal generated by
\emph{isotopy} (not displayed; we impose all possible plane isotopies), \emph{circle} and \emph{bitri evaluation}, and the \emph{H=I relation}:
\begin{gather*}
\begin{tikzpicture}[anchorbase,scale=1]
\draw[usual] (0,0)to[out=90,in=180] (0.5,0.5)to[out=0,in=90] (1,0);
\draw[usual] (0,0)to[out=270,in=180] (0.5,-0.5)to[out=0,in=270] (1,0);
\end{tikzpicture}
=\qnum[3]
=\qpar^{2}+1+\qpar^{-2}
,\quad
\begin{tikzpicture}[anchorbase,scale=1]
\draw[usual] (0,0)to[out=90,in=180] (0.5,0.5);
\draw[usual] (1,0)to[out=90,in=0] (0.5,0.5);
\draw[usual] (0.5,0)to (0.5,0.5);
\draw[usual] (1,-0.5) to (1,0);
\draw[usual] (0,0)to[out=270,in=180] (0.25,-0.5)to[out=0,in=270] (0.5,0);
\end{tikzpicture}
=0
,\\
\begin{tikzpicture}[anchorbase,scale=1]
\draw[usual] (0,0)to (0.33,0.5)to (0,1);
\draw[usual] (1,0)to (0.66,0.5)to (1,1);
\draw[usual] (0.33,0.5)to (0.66,0.5);
\end{tikzpicture}
=
\begin{tikzpicture}[anchorbase,scale=1]
\draw[usual] (0,0)to (0.5,0.33)to (1,0);
\draw[usual] (0,1)to (0.5,0.66)to (1,1);
\draw[usual] (0.5,0.33)to (0.5,0.66);
\end{tikzpicture}
+1/(\qpar^{2}+\qpar^{-2})\cdot
\begin{tikzpicture}[anchorbase,scale=1]
\draw[usual] (0,0)to[out=45,in=315] (0,1);
\draw[usual] (1,0)to[out=135,in=225] (1,1);
\end{tikzpicture}
-1/(\qpar^{2}+\qpar^{-2})\cdot
\begin{tikzpicture}[anchorbase,scale=1]
\draw[usual] (0,0)to[out=45,in=135] (1,0);
\draw[usual] (0,1)to[out=315,in=225] (1,1);
\end{tikzpicture}
\,.
\end{gather*}
We call $\web[{\son[3]}]$ the \emph{$\son[3]$ web category} and 
its morphism \emph{$\son[3]$ webs}.
\end{Definition}

\begin{Remark}\label{R:RepsSO3NotDefined}
Note that we do not define the $\son[3]$ web 
category for $\qpar^{2}+\qpar^{-2}=0$. In particular, when talking 
about this category we will always assume that 
$\qpar^{2}+\qpar^{-2}\neq 0$.
\end{Remark}

\begin{Example}\label{E:RepsSO3Relation}
The H=I relation can be used to systematically reduce faces of $\son[3]$ webs 
in their complexity. For example,
\begin{gather*}
\scalebox{0.95}{$\begin{tikzpicture}[anchorbase,scale=1]
\draw[usual] (0,0)to[out=90,in=180] (0.5,0.5);
\draw[usual] (1,0)to[out=90,in=0] (0.5,0.5);
\draw[usual] (0.5,0)to (0.5,0.5);
\draw[usual] (0.5,-0.5)to[out=180,in=270] (0,0);
\draw[usual] (0.5,-0.5)to[out=0,in=270] (1,0);
\draw[usual] (0.5,-0.5)to (0.5,0);
\draw[very thick,densely dotted,blue] (0.3,-0.6) rectangle (0.7,0.6);
\end{tikzpicture}
=
\underbrace{\begin{tikzpicture}[anchorbase,scale=1]
\draw[usual] (0,0)to (0.33,0.5)to (0,1)to[out=180,in=180] (0,0);
\draw[usual] (1,0)to (0.66,0.5)to (1,1)to[out=0,in=0] (1,0);
\draw[usual] (0.33,0.5)to (0.66,0.5);
\end{tikzpicture}}_{=0}
-1/(\qpar^{2}+\qpar^{-2})\cdot
\underbrace{\begin{tikzpicture}[anchorbase,scale=1]
\draw[usual] (0,0)to[out=90,in=180] (0.5,0.5)to[out=0,in=90] (1,0);
\draw[usual] (0,0)to[out=270,in=180] (0.5,-0.5)to[out=0,in=270] (1,0);
\draw[usual] (1.25,0)to[out=90,in=180] (1.75,0.5)to[out=0,in=90] (2.25,0);
\draw[usual] (1.25,0)to[out=270,in=180] (1.75,-0.5)to[out=0,in=270] (2.25,0);
\end{tikzpicture}}_{=\qnum[3]^{2}}
+1/(\qpar^{2}+\qpar^{-2})\cdot
\underbrace{\begin{tikzpicture}[anchorbase,scale=1]
\draw[usual] (0,0)to[out=90,in=180] (0.5,0.5)to[out=0,in=90] (1,0);
\draw[usual] (0,0)to[out=270,in=180] (0.5,-0.5)to[out=0,in=270] (1,0);
\end{tikzpicture}}_{=\qnum[3]}
=-\qnum[3].$}
\end{gather*}
In the first picture we highlighted 
an I that we then replaced by H and error terms.
In a similar fashion one can get relations for other faces as well.
\end{Example}

\ochanged{A \emph{higher valent vertices}, exemplified, is:
\begin{gather*}
\emph{4 valent}\colon
\begin{tikzpicture}[anchorbase,scale=1]
\draw[usual] (0,0)to (1,1);
\draw[usual] (1,0)to (0,1);
\draw[tomato,fill=tomato] (0.5,0.5) circle (0.1cm);
\end{tikzpicture}
,\quad
\emph{7 valent}\colon
\begin{tikzpicture}[anchorbase,scale=1]
\draw[usual] (0,0)to (0.43,-0.9);
\draw[usual] (0,0)to (0.97,-0.22);
\draw[usual] (0,0)to (0.78,0.62);
\draw[usual] (0,0)to (0,1);
\draw[usual] (0,0)to (-0.78,0.62);
\draw[usual] (0,0)to (-0.97,-0.22);
\draw[usual] (0,0)to (-0.43,-0.9);
\draw[tomato,fill=tomato] (0,0) circle (0.1cm);
\end{tikzpicture}
\,.
\end{gather*}
The dot is a visual aid only.
Diagrams that are allowed to have these additional vertices are embedded graphs with specified bottom and top boundary.
An edge of such a graph is called \emph{inner} if it does not touch the 
boundary. The \emph{contraction} operation is
\begin{gather*}
\begin{tikzpicture}[anchorbase,scale=1]
\draw[usual] (0,0)to (0.33,0.5)to (0,1);
\draw[usual] (1,0)to (0.66,0.5)to (1,1);
\draw[usual,dotted] (0.33,0.5)to (0.66,0.5);
\end{tikzpicture}
\mapsto
\begin{tikzpicture}[anchorbase,scale=1]
\draw[usual] (0,0)to (1,1);
\draw[usual] (1,0)to (0,1);
\draw[tomato,fill=tomato] (0.5,0.5) circle (0.1cm);
\end{tikzpicture}
\reflectbox{$\mapsto$}
\begin{tikzpicture}[anchorbase,scale=1]
\draw[usual] (0,0)to (0.5,0.33)to (1,0);
\draw[usual] (0,1)to (0.5,0.66)to (1,1);
\draw[usual,dotted] (0.5,0.33)to (0.5,0.66);
\end{tikzpicture}
,
\end{gather*}
where the dotted edge is contracted.}

\begin{Definition}\label{D:RepsSO3Basis}
Let $k+l$ be the number of boundary points of $\son$ web $u$.
We say $u$ is a \emph{partition} (of the set $\{1,\dots,k+l\}$) if:
\begin{enumerate}

\item $u$ is one $\son$ web, i.e. not a nontrivial $\C$-linear combination of such diagrams.

\item \ochanged{$u$ has no internal faces.}

\item \ochanged{After a finite number of contractions, $u$ is a graph without inner edges. (Here we see $u$ as a trivalent graph and then apply contraction.)}

\end{enumerate}
Let $\sonbasis{k}{l}$ be the set of all partition $\son$ web diagrams
with $k$ bottom and $l$ top boundary points.
\end{Definition}

\begin{Lemma}\label{L:RepsSO3Basis}
The set $\sonbasis{k}{l}$ is a $\C$-basis of $\Hom_{\wweb[{\son}]}(\genob^{\hcirc k},\genob^{\hcirc l})$.
\end{Lemma}

\begin{proof}
\textit{Spanning.} As exemplified in \autoref{E:RepsSO3Relation}, the H=I relation implies that we can assume that $u$ has no internal faces. Indeed, the faces marked with a bullet in
\begin{gather*}
\begin{tikzpicture}[anchorbase,scale=1]
\draw[usual] (0,0)to (0.33,0.5)to (0,1);
\draw[usual] (1,0)to (0.66,0.5)to (1,1);
\draw[usual] (0.33,0.5)to (0.66,0.5);
\node[circle,draw=black,fill=tomato,inner sep=0pt,minimum size=5pt] at (0.5,0.9) {};
\end{tikzpicture}
=
\begin{tikzpicture}[anchorbase,scale=1]
\draw[usual] (0,0)to (0.5,0.33)to (1,0);
\draw[usual] (0,1)to (0.5,0.66)to (1,1);
\draw[usual] (0.5,0.33)to (0.5,0.66);
\node[circle,draw=black,fill=tomato,inner sep=0pt,minimum size=5pt] at (0.5,0.9) {};
\end{tikzpicture}
\,,\quad
\begin{tikzpicture}[anchorbase,scale=1]
\draw[usual] (0,0)to (0.5,0.33)to (1,0);
\draw[usual] (0,1)to (0.5,0.66)to (1,1);
\draw[usual] (0.5,0.33)to (0.5,0.66);
\node[circle,draw=black,fill=tomato,inner sep=0pt,minimum size=5pt] at (0.15,0.5) {};
\end{tikzpicture}
=
\begin{tikzpicture}[anchorbase,scale=1]
\draw[usual] (0,0)to (0.33,0.5)to (0,1);
\draw[usual] (1,0)to (0.66,0.5)to (1,1);
\draw[usual] (0.33,0.5)to (0.66,0.5);
\node[circle,draw=black,fill=tomato,inner sep=0pt,minimum size=5pt] at (0.1,0.5) {};
\end{tikzpicture}
\,,
\end{gather*}
will have fewer edges on the \tchanged{right-hand sides when compared to the left-hand sides. We can repeat this operation until some internal face is a monogon and the bitri evaluation applies.}
Moreover, the two error terms in the H=I relation are simpler $\son$ webs since the number of vertices is smaller than for the other two $\son$ webs.
In other words, internal faces can be \tchanged{removed} recursively. 
Finally, the H=I relation let us get rid of inner edges, which shows that 
$\sonbasis{k}{l}$ spans.

\textit{Linear independence.} There is a bijection from $\sonbasis{k}{l}$ 
to the set of all planar partitions of the set $\{1,\dots,k+l\}$ where every block has at least two parts given by associating a partition 
to a partition $\son$ web diagram by interpreting the connected components of the web as blocks of the partition. Let $pp(k,l)$ be the number of such partitions.
Since $\sonbasis{k}{l}$ spans 
$\Hom_{\wweb[{\son}]}(\genob^{\otimes k},\genob^{\otimes l})$, we get
$\dim_{\C}\Hom_{\wweb[{\son}]}(\genob^{\otimes k},\genob^{\otimes l})\leq pp(k,l)$, while pivotality 
and \cite[Lemma 4.1]{FeKr-chomatic-category} imply that 
$pp(k,l)\leq\dim_{\C}\Hom_{\wweb[{\son}]}(\genob^{\otimes k},\genob^{\otimes l})$. Hence, linear independence follows.
\end{proof}

\begin{Remark}\label{R:RepsSO3Basis}
\ochanged{The numbers $pp(k,l)$ are well-known in combinatorics. Without loss of generality we can consider $pp(k,l)$ for $l=0$ and one gets
\begin{gather*}
\{1,0,1,1,3,6,15,36,91,232,603\},\quad pp(k,0)\text{ for }k=0,\dots,10.
\end{gather*}
The sequence is \cite[A005043]{oeis}.}
\end{Remark}

Let us denote symmetric 
$\sln$ webs in the sense of \cite{RoTu-symmetric-howe} by using labeled (and colored) edges, for example,
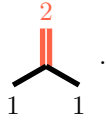
\begin{gather*}
\begin{tikzpicture}[anchorbase,scale=1]
\draw[symedge] (0,0.5) node[above,yshift=-1pt]{$2$}to (0,0);
\draw[usual] (-0.433013,-0.25)node[below]{$1$}to (0,0);
\draw[usual] (0.433013,-0.25)node[below]{$1$}to (0,0);
\end{tikzpicture}
\,.
\end{gather*}
The edge labeled $1$ are uncolored.
Let $\symweb$ denote the associated $\C$-linear pivotal category.

\begin{Lemma}\label{L:RepsSO3SL2}
Assume $\qpar\in\C\setminus\{0\}$ is not a nontrivial 
root of unity.
There is a faithful $\C$-linear pivotal functor
$\functorstuff{I}\colon\web[{\son}]\to\symweb$ determined by
\begin{gather*}
\begin{tikzpicture}[anchorbase,scale=1]
\draw[usual] (0,0)to[out=90,in=180] (0.5,0.5)to[out=0,in=90] (1,0);
\end{tikzpicture}
\mapsto
\begin{tikzpicture}[anchorbase,scale=1]
\draw[symedge] (0,0)node[below]{$2$}to[out=90,in=180] (0.5,0.5)node[above]{\phantom{2}}to[out=0,in=90] (1,0)node[below]{$2$};
\end{tikzpicture}
\,,\quad
\begin{tikzpicture}[anchorbase,scale=1]
\draw[usual] (0,0)to[out=270,in=180] (0.5,-0.5)to[out=0,in=270] (1,0);
\end{tikzpicture}
\mapsto
\begin{tikzpicture}[anchorbase,scale=1]
\draw[symedge] (0,0)node[above,yshift=-1pt]{$2$}to[out=270,in=180] (0.5,-0.5)node[below]{\phantom{2}}to[out=0,in=270] (1,0)node[above,yshift=-1pt]{$2$};
\end{tikzpicture}
\,,\\
\begin{tikzpicture}[anchorbase,scale=1]
\draw[usual] (0,0)to[out=90,in=180] (0.5,0.5);
\draw[usual] (1,0)to[out=90,in=0] (0.5,0.5);
\draw[usual] (0.5,0)to (0.5,0.5);
\end{tikzpicture}
\mapsto
\tfrac{1}{s}\cdot
\begin{tikzpicture}[anchorbase,scale=1]
\draw[symedge] (0,0)node[below]{$2$}to (0,0.25);
\draw[symedge] (0.5,0)node[below]{$2$}to (0.5,0.25);
\draw[symedge] (1,0)node[below]{$2$}to (1,0.25);
\draw[usual] (0,0.25)to[out=45,in=135](0.5,0.25)to[out=45,in=135]
(1,0.25)to[out=45,in=315](0.85,0.65)to[out=135,in=45](0.15,0.65)node[above]{\phantom{1}}to[out=225,in=135](0,0.25);
\end{tikzpicture}
\,,\quad
\begin{tikzpicture}[anchorbase,scale=1]
\draw[usual] (0.5,-0.5)to[out=180,in=270] (0,0);
\draw[usual] (0.5,-0.5)to[out=0,in=270] (1,0);
\draw[usual] (0.5,-0.5)to (0.5,0);
\end{tikzpicture}
\mapsto
\tfrac{1}{s}\cdot
\begin{tikzpicture}[anchorbase,scale=1,yscale=-1]
\draw[symedge] (0,0)node[above,yshift=-1pt]{$2$}to (0,0.25);
\draw[symedge] (0.5,0)node[above,yshift=-1pt]{$2$}to (0.5,0.25);
\draw[symedge] (1,0)node[above,yshift=-1pt]{$2$}to (1,0.25);
\draw[usual] (0,0.25)to[out=45,in=135](0.5,0.25)to[out=45,in=135]
(1,0.25)to[out=45,in=315](0.85,0.65)to[out=135,in=45](0.15,0.65)node[above]{\phantom{1}}to[out=225,in=135](0,0.25);
\end{tikzpicture}
\,,
\end{gather*}
where we choose a square root $s=\big((\qpar^{2}+\qpar^{-2})\qnum[2]^{2}\big)^{1/2}$ of 
$(\qpar^{2}+\qpar^{-2})\qnum[2]^{2}$.
\end{Lemma}

\begin{proof}
A direct verification shows that the above defines a $\C$-linear pivotal functor.
There are only two things to note here. Firstly, the scaling which comes from the comparison of the relations
\begin{gather*}
\begin{tikzpicture}[anchorbase,scale=1]
\draw[usual] (0,0)to[out=90,in=180] (0.5,0.5);
\draw[usual] (1,0)to[out=90,in=0] (0.5,0.5);
\draw[usual] (0.5,0)to (0.5,0.5);
\draw[usual] (0.5,-0.5)to[out=180,in=270] (0,0);
\draw[usual] (0.5,-0.5)to[out=0,in=270] (1,0);
\draw[usual] (0.5,-0.5)to (0.5,0);
\end{tikzpicture}
=-\qnum[3]
\text{ and }
\begin{tikzpicture}[anchorbase,scale=1]
\draw[symedge] (0,-0.25)to (0,0)node[left]{$2$}to (0,0.25);
\draw[symedge] (0.5,-0.25)to (0.5,0)to (0.5,0.25);
\draw[symedge] (1,-0.25)to (1,0)to (1,0.25);
\draw[usual] (0,0.25)to[out=45,in=135](0.5,0.25)to[out=45,in=135]
(1,0.25)to[out=45,in=315](0.85,0.65)to[out=135,in=45](0.15,0.65)node[above]{\phantom{1}}to[out=225,in=135](0,0.25);
\draw[usual] (0,-0.25)to[out=-45,in=-135](0.5,-0.25)to[out=-45,in=-135]
(1,-0.25)to[out=-45,in=-315](0.85,-0.65)to[out=-135,in=-45](0.15,-0.65)node[below]{\phantom{1}}to[out=-225,in=-135](0,-0.25);
\end{tikzpicture}
=-(\qpar^{2}+\qpar^{-2})\qnum[2]^{2}\qnum[3].
\end{gather*}
Second, to verify the defining relations hold in the image of $\functorstuff{I}$ is an easy calculation.

That $\functorstuff{I}$ is an embedding can be checked by using 
\autoref{L:RepsSO3Basis} and the 
faithful representation $\Gamma_{sym}$ of $\symweb$ on symmetric powers obtained from the functor used in the proof of
\cite[Theorem 1.10]{RoTu-symmetric-howe}. The only 
thing the reader needs to know to verify this is the following. Fix 
the basis $\{v_{1},v_{2}\}$ of $\C^{2}$. The basis elements of 
$\mathrm{Sym}^{2}\C^{2}$ are chosen to be $\{v_{1}v_{1},v_{1}v_{2}=\qpar^{-1}\cdot v_{2}v_{1},v_{2}v_{2}\}$.
Then
\begin{gather*}
\begin{tikzpicture}[anchorbase,scale=1]
\draw[usual] (0,0)node[below]{$1$}to[out=90,in=180] (0.5,0.5)node[above]{\phantom{1}}to[out=0,in=90] (1,0)node[below]{$1$};
\end{tikzpicture}
\mapsto
\left\{
\begin{gathered}
v_{i}\hcirc v_{i}\mapsto 0,
\\
v_{1}\hcirc v_{2}\mapsto -\qpar,
v_{2}\hcirc v_{1}\mapsto 1,
\end{gathered}
\right.
\\
\begin{tikzpicture}[anchorbase,scale=1,yscale=-1]
\draw[usual] (0,0)node[above,yshift=-1pt]{$1$}to[out=90,in=180] (0.5,0.5)node[below]{\phantom{1}}to[out=0,in=90] (1,0)node[above,yshift=-1pt]{$1$};
\end{tikzpicture}
\mapsto
(1\mapsto v_{1}\hcirc v_{2}-\qpar^{-1}\cdot v_{2}\hcirc v_{1})
,\\
\begin{tikzpicture}[anchorbase,scale=1]
\draw[symedge] (0,0.5) node[above,yshift=-1pt]{$2$}to (0,0);
\draw[usual] (-0.433013,-0.25)node[below]{$1$}to (0,0);
\draw[usual] (0.433013,-0.25)node[below]{$1$}to (0,0);
\end{tikzpicture}
\mapsto
(v_{i}\hcirc v_{j}\mapsto v_{i}v_{j})
,\quad
\begin{tikzpicture}[anchorbase,scale=1,yscale=-1]
\draw[symedge] (0,0.5) node[below]{$2$}to (0,0);
\draw[usual] (-0.433013,-0.25)node[above,yshift=-1pt]{$1$}to (0,0);
\draw[usual] (0.433013,-0.25)node[above,yshift=-1pt]{$1$}to (0,0);
\end{tikzpicture}
\mapsto
\left\{
\begin{gathered}
v_{i}v_{i}\mapsto\qnum\cdot v_{i}\hcirc v_{i},
\\
v_{1}v_{2}\mapsto\qpar^{-1}v_{1}\hcirc v_{2}+v_{2}\hcirc v_{1}\ochanged{,}
\end{gathered}
\right.
\end{gather*}
under $\Gamma_{sym}$, while the $2$ labeled caps and cups are defined by explosion, see \cite[Definition 2.18]{RoTu-symmetric-howe}.
\end{proof}

\begin{Lemma}\label{L:RepsSO3Semisimple}
We have the following.
\begin{enumerate}

\item The simple objects of $\web[{\son[3]}]$ are in one-to-one correspondence 
with $\N$.

\item $\web[{\son[3]}]$ is semisimple
if only if $\qpar\in\C\setminus\{0\}$ is not a nontrivial 
root of unity.

\item $\web[{\son[3]}]$ is \changed{an essentially fusion} category if only if $\qpar\in\C\setminus\{0\}$ is not a nontrivial 
root of unity.

\end{enumerate}
\end{Lemma}

\begin{proof}
This follows from \autoref{L:RepsSL2Semisimple} 
and the fact that $\son$ webs can be constructed as the full subcategory 
of $\sln$ webs $\hcirc$-generated by the diagrammatic 
analog of $\mathrm{Sym}^{2}\C^{2}$, see \autoref{L:RepsSO3SL2}.
\end{proof}

For $\son[3]$ webs the crossing formulas are:
\begin{gather}\label{Eq:RepsSO3Braiding}
\begin{gathered}
\begin{tikzpicture}[anchorbase,scale=1]
\draw[usual] (1,0)to (0,1);
\draw[usual,crossline] (0,0)to (1,1);
\end{tikzpicture}
=
(\qpar^{2}-1)\cdot
\begin{tikzpicture}[anchorbase,scale=1]
\draw[usual] (0,0)to[out=45,in=315] (0,1);
\draw[usual] (1,0)to[out=135,in=225] (1,1);
\end{tikzpicture}
+\qpar^{-2}\cdot
\begin{tikzpicture}[anchorbase,scale=1]
\draw[usual] (0,0)to[out=45,in=135] (1,0);
\draw[usual] (0,1)to[out=315,in=225] (1,1);
\end{tikzpicture}
+(\qpar^{2}+\qpar^{-2})\cdot
\begin{tikzpicture}[anchorbase,scale=1]
\draw[usual] (0,0)to (0.33,0.5)to (0,1);
\draw[usual] (1,0)to (0.66,0.5)to (1,1);
\draw[usual] (0.33,0.5)to (0.66,0.5);
\end{tikzpicture}
\,,
\\
\begin{tikzpicture}[anchorbase,scale=1]
\draw[usual] (0,0)to (1,1);
\draw[usual,crossline] (1,0)to (0,1);
\end{tikzpicture}
=
(\qpar^{-2}-1)\cdot
\begin{tikzpicture}[anchorbase,scale=1]
\draw[usual] (0,0)to[out=45,in=315] (0,1);
\draw[usual] (1,0)to[out=135,in=225] (1,1);
\end{tikzpicture}
+\qpar^{2}\cdot
\begin{tikzpicture}[anchorbase,scale=1]
\draw[usual] (0,0)to[out=45,in=135] (1,0);
\draw[usual] (0,1)to[out=315,in=225] (1,1);
\end{tikzpicture}
+(\qpar^{2}+\qpar^{-2})\cdot
\begin{tikzpicture}[anchorbase,scale=1]
\draw[usual] (0,0)to (0.33,0.5)to (0,1);
\draw[usual] (1,0)to (0.66,0.5)to (1,1);
\draw[usual] (0.33,0.5)to (0.66,0.5);
\end{tikzpicture}
\,.
\end{gathered}
\end{gather}

Now all of \autoref{L:RepsSL2Braiding} 
(with \autoref{Eq:RepsSO3Braiding}) and \autoref{N:RepsSL2Braiding} have the 
evident $\son[3]$ analog (their formulation is omitted) and we will use these analogs freely.
\tchanged{In particular, $\web[{\son[3]}]$ is a braided category.}


\subsection{The main theorem in the SO3 case}\label{SS:RepsSO3Main}


Recall that we introduced our notation for tensors in \autoref{SS:RepsGL2Main}. We will use the same conventions now.

\begin{Theorem}\label{T:RepsSO3Main}
Assume $\qpar\in\C\setminus\{0\}$ is not a nontrivial 
root of unity.
\begin{enumerate}

\item Let $n\geq 3$. For every pair $\bN=(\mat,\ten)\in\matgln\times\tensors{n}{n}{n}$ with $\trace(\mat^{T}\mat^{-1})=\qnum[3]$, 
$\trace\big(\ten(\mat(cup)\hcirc id)\big)=0$, 
and 
$(id\hcirc\ten_{l})\vcirc(\ten^{l}\hcirc id)=\ten^{l}\ten_{l}+1/(\qpar^{2}+\qpar^{-2})\cdot id-1/(\qpar^{2}+\qpar^{-2})\cdot\mat(cap)\mat(cup)$ there exists 
a simple transitive fiber 2-representation $\mrep_{\mat,\ten}^{n}$ of $\web[{\son[3]}]$ constructed \changed{similarly} to the proof of \autoref{L:RepsSL2Matrix}.
\textbf{(Existence)}

\item We have $\mrep_{(\mat,\ten)}\equirep\mrep_{(\mat[M],\ten[U])}$	
if and only if $\mat\congruence\mat[M]$ 
and $\ten\congruence\mat[U]$.  
\textbf{(Non-redundant)}

\item Every simple transitive fiber 2-representation of $\web[{\son[3]}]$
is of the form $\mrep_{\bN}$, and every simple transitive rank one 2-representation of $\web[{\son[3]}]$ arises in this way.
\textbf{(Complete)}

\end{enumerate}
Moreover, there are infinitely many nonequivalent 
simple transitive rank one 2-representations of $\web[{\son}]$.
\end{Theorem}

We also show that for $n=3$ there is only one possible solution.

Essentially all we said at the end of \autoref{SS:RepsSL2Main} 
(before the remarks) 
applies for $\son$ webs as well. In particular, 
we leave the analog of \autoref{P:RepsSL2Main} to the reader, and we 
will only focus on the crucial difference compared to the other two cases: the 
appearance of (honest) trilinear forms. This might make a \ochanged{``huge''} difference, see \autoref{S:Complexity} 
for a more detailed discussion.

\begin{Remark}\label{R:RepsSO3HistRemark}
The category $\web[{\son[3]}]$ was discovered in the early days of quantum topology, see \cite{Ya-invariant-graphs} for the potentially earliest reference. 
In that paper it is effectively shown that $\web[{\son[3]}]$ gives a diagrammatic 
description of $\son$-representations (this can be pieced together 
by comparing \autoref{L:RepsSO3SL2} and the MathSciNet review of 
\cite{Ya-invariant-graphs}). As far as we know, $\web[{\son[3]}]$ is the oldest diagram category that truly deserves the name web category.
Its importance stems from its connection to, for example, 
the chromatic polynomial and the four color theorem in graph theory.
This connection originates in \cite{TeLi-the-tl-paper}, see \cite[Introduction]{FeKr-tutte-tl-algebra} 
for a list of early appearances of this relation.
\end{Remark}

\begin{Remark}\label{R:RepsSO3HistRemarkTwo}
In contrast to \autoref{T:RepsSL2Main}, a generalization 
of \autoref{T:RepsSO3Main} beyond rank one appears to be difficult.
See however \cite{EvPu-so3-classification} for a related classification.
\end{Remark}

\begin{Remark}\label{R:RepsSO3HistRemarkThree}
\autoref{T:RepsSO3Main} seems very different than
\cite[Theorems 1.1 and 1.2]{Mr-quantum-so3}.
\end{Remark}


\subsection{Proof of \autoref{T:RepsSO3Main}}\label{SS:RepsSO3Proof}


A tuple $\bN=(\mat,\ten)\in\matgln\times\tensors{n}{n}{n}$ 
as in \autoref{T:RepsSO3Main} is called an \textit{$\son$ tuple.}

\begin{Lemma}\label{L:RepsSO3Matrix}
For $m\in\Z_{\geq 2}$ let 
$\mat\in\matgln[m]$ be a matrix satisfying $\trace(\mat^{T}\mat^{-1})=-\qnum$. 
Then there exists an associated $\son$ tuple with $n=m+1$.
\end{Lemma}

\begin{proof}
Recall \ochanged{from \cite[Proof of Theorem 1.10]{RoTu-symmetric-howe}} that $\symweb$ is monoidally equivalent to 
$\web$ upon taking additive idempotent closures, and the equivalence is given 
by an explicit monoidal functor $\functorstuff{F}$. \changed{In a bit more detail}, the object $k$ in 
$\symweb$ corresponds to the $k$th Jones--Wenzl projector in $\web$. In any case, 
we get a monoidal equivalence $\functorstuff{F}\colon\addc{\web}\to\add{\symweb}\cong_{\hcirc}\addc{\symweb}$
between the additive idempotent closure of $\web$ and the additive closure of $\symweb$. We identify the two categories using $\functorstuff{F}$.

Recall $\functorstuff{I}$ from \autoref{L:RepsSO3SL2},
and consider the following commutative diagram:
\begin{gather*}
\begin{tikzpicture}[anchorbase,->,>=stealth',shorten >=1pt,auto,node distance=3cm,
thick,main node/.style={font=\sffamily\Large\bfseries}]
\node[main node] (1) {$\web[{\son}]$};
\node[main node] (2) [right of=1,xshift=0.5cm] {$\addc{\web}$};
\node[main node] (3) [below of=2,yshift=1.5cm] {$\web$};
\node[main node] (4) [right of=2,xshift=0.5cm] {$\vect$};
\node[main node] (5) [right of=3,xshift=0.5cm] {$\vect$.};
\path[every node/.style={font=\sffamily\small}]
(1) edge[right hook-stealth] node[above]{$\functorstuff{I}$} (2)
(2) edge[orchid,dotted] node[above]{$\exists!\tilde{\mrep}_{\mat}^{n}$} (4)
(3) edge[right hook-stealth] node[left]{$\text{incl.}$} (2)
(3) edge node[below]{$\mrep_{\mat}^{n}$} (5)
(4) edge[double equal sign distance,-] (5);
\end{tikzpicture}
\end{gather*}
The existence of $\tilde{\mrep}_{\mat}^{n}$ follows 
from the usual \changed{Yoga of additive
and idempotent closures}. Thus, we get a 2-representation $\tilde{\mrep}_{\mat}^{n}\vcirc\functorstuff{I}$ of $\web[{\son}]$.

Note that all needed functors are given explicitly. Tracking back their 
definitions and a bit of calculation gives the desired matrices and tensors.
\end{proof}

\begin{Lemma}\label{L:RepsSO3MatrixTwo}
For every $n\in\Z_{\geq 3}$ there exists some
$\son$ tuple.
For $n\in\{1,2\}$ there exists no such tuples.
\end{Lemma}

\begin{proof}
By \autoref{L:RepsSO3Matrix} and the corresponding statement for $\web$ as in \autoref{L:RepsSL2MatrixTwo}, we get the existence. The case $n=1$ is ruled out 
as in \autoref{R:RepsSL2Reps}, while $n=2$ can be ruled out as in
\autoref{L:RepsSL2nTwo}.
\end{proof}

\begin{Lemma}\label{L:RepsSO3Collect}
The statements \autoref{L:RepsSL2MatrixEquivalence} to 
\autoref{L:RepsSL2MatrixEquiClasses}
hold mutatis mutandis for $\web[{\son}]$ as well.
\end{Lemma}

\begin{proof}
\tchanged{Let us go through the lemmas one-by-one and mention what needs to be changed:
\begin{enumerate}
\item For \autoref{L:RepsSL2MatrixEquivalence} we first recall that \autoref{L:RepsSO3MatrixTwo} shows that for $n=3$ the only possible $\son$ tuple is the one coming from the standard choice (given by $\son$ acting on its defining representation), while there are no solutions 
for $n<3$. Moreover, for $n=3$ the lemma follows by using 
\autoref{L:RepsSO3SL2} and then tracking the image of the basis from \autoref{L:RepsSO3Basis} under quantum symmetric Howe duality.The general case follows by copying this for a higher dimensional target space.
\item \autoref{L:RepsSL2MatrixAll} works in the same way: one chooses a basis and orders the images of the generators in corresponding matrices or tensors.
\item In \autoref{L:RepsSL2MatrixSimple} one replaces 
the reference to \cite[Chapter XII]{Tu-qgroups-3mfds} 
with \cite[Lemma 3.4]{Mr-quantum-so3}. Indeed, the proof of \cite[Lemma 3.4]{Mr-quantum-so3} can be copied as it only relies on the fusion rules of $\son$. We get the desired unique functor, up to scaling, as all generators exists uniquely, up to scaling, as maps and all relations are satisfied, for example, 
\cite[Lemma 3.4.(2f)]{Mr-quantum-so3} is the H=I relation.
\item \autoref{L:RepsSL2Same} follows as before from the previous two points.
\item Ditto, \autoref{L:RepsSL2MatrixEquiClasses} follows as before from the third point above.
\end{enumerate}
Details are omitted.}
\end{proof}

We have a complete solution for matrix congruence, see \autoref{SS:RepsSL2Proof},
which is the same as equivalence of bilinear forms by the classical fact that two matrices are congruent if and only if they represent the same bilinear form up to change-of-basis. 

\begin{Example}\label{E:RepSO3List}
For $n=3$ \autoref{L:RepsSL2List} lists all possible solutions of $\trace(\mat^{T}\mat^{-1})=-\qnum$ up to $\congruence$. 
The solutions of $\trace(\mat^{T}\mat^{-1})=\qnum[3]$ up to $\congruence$ 
are similar. That is, for $\qpar$ generic enough the only possible solution is
\begin{gather*}
\mat=\mat[G]_{1}\oplus\mat[H]_{1}(\qpar^{2})=
\begin{psmallmatrix}
1 & 0 & 0
\\
0 & 0 & 1
\\
0 & \qpar^{2} & 0
\end{psmallmatrix},
\end{gather*}
up to $\qpar\leftrightarrow\qpar^{-1}$. This is the standard solution up to permutation.
\end{Example}

Thus, we only need to worry about trilinear forms. The easiest case for us 
are \emph{ternary trilinear forms}, often called $(3,3,3)$ trilinear forms,
where $n=3$. In the notation above this is the case displayed in \autoref{Eq:RepsSO3Trilinear}. For $\bx[1]=(1,1,1)$, we denote the appearing $3$-by-$3$ matrices by $\ten_{x}(\bx[1])$, $\ten_{y}(\bx[1])$ 
and $\ten_{z}(\bx[1])$ in order from front to back.

Take now such a form $\ten$ and write it as $\ten=\sum_{h,i,j}t_{hij}\cdot x_{h}y_{i}z_{j}$, using variables. We let $\ten_{x}(\bx)=(\sum_{h}t_{hij}\cdot x_{h})_{ij}$ 
for $\bx=(x_{1},x_{2},x_{3})$, and similarly 
$\ten_{y}(\bx[y])$ and $\ten_{z}(\bx[z])$. The determinant formula $\det\big(\ten_{x}(\bx)\big)=0$ 
is a ternary cubic that we denote by $\ten_{x}$. We also have $\ten_{y}$ and $\ten_{z}$ by using the corresponding determinant formulas.
Finally, evaluation at $\bx[a]\in\C^{3}$ gives $\ten_{x}(\bx[a])$. This is a complex matrix, so we can let $t_{x}\in\N\cup\{\infty\}$ be the number of matrices $\ten_{x}(\bx[a])$ with $\rank_{\C}\ten_{x}(\bx[a])=1$ (this number 
can be infinite). Similarly for $t_{y}$ and $t_{z}$.

\begin{Lemma}\label{L:RepsSO3TrilinearPreLim}
Any ternary cubic is projectively equivalent to one
of the following:
\begin{gather*}
1:x^{3}=0,\quad
2:x^{2}y=0,\quad
3:xy(x-y)=0,\quad
4:xyz=0,\\
5:z(x^{2}+yz)=0,\quad
6:x(x^{2}+yz)=0,\quad
7:x^{3}-y^{2}z=0,\quad
8:x^{3}+y^{3}-xyz=0,
\end{gather*}
as well as $9:$ an elliptic cubic and $10:$ a zero cubic.
\end{Lemma}

\begin{proof}
Well-known, see \cite{ThCh-trilinear}. A more modern 
and detailed account can be found in many works, see for example
\cite[Table 1]{LaTe-ranks-tensors}.
\end{proof}

One has a complete classification of ternary trilinear forms:

\begin{Lemma}\label{L:RepsSO3Trilinear}
We have the following.
\begin{enumerate}

\item We have $\ten\congruence\mat[U]$ if and only if ($(t_{x},t_{y},t_{z})$ 
is equal to $(u_{x},u_{y},u_{z})$ in some order, and $(\ten_{x},\ten_{y},\ten_{z})$
is projective equivalent to $(\mat[U]_{x},\mat[U]_{y},\mat[U]_{z})$ in the same order.)

\item The only possible triples $(t_{x},t_{y},t_{z})$, up to reordering, are 
listed in the table in \autoref{Eq:RepsSO3Table} below. The only possible ternary cubics, up to projective equivalence, are listed in the same table.

\end{enumerate}
\begin{gather}\label{Eq:RepsSO3Table}
\begin{gathered}
\scalebox{0.8}{\begin{tabular}{C||C|C|C|C|C|C|C|C|C}
& (0,1,0) & (1,0,1) & (1,1,1) & (1,2,1) & (2,1,2) & (2,2,2) & (3,3,3) & (\infty,1,\infty) & (\infty,2,\infty) \\
\hline
\hline
1 &  &  & 1 &  &  &  &  & 10 &  \\
\hline
2 &  &  & 2,5 & 1 & 3 & 2 &  &  & 10 \\
\hline
3 &  &  &  & 2 &  &  &  &  &  \\
\hline
4 &  &  &  4 & 6 &  & 4 & 4 & & \\
\hline
5 & 7 & 3 &  &  &  &  &  &  &  \\
\hline
6 & 8 & 4 &  &  &  &  &  &  &  \\
\end{tabular}}
\\
(0,0,0)\colon\text{nonzero $\Leftrightarrow$ all cubics are of the same projective type.}
\end{gathered}
\end{gather}
\autoref{Eq:RepsSO3Table} is to be read as follows. The projective cases of 
two of $(\ten_{x},\ten_{y},\ten_{z})$ need to 
agree, up to order, and the 
list in the first column is the class for $\ten_{x}$ and $\ten_{z}$. The class of $\ten_{y}$ is then listed in the table, with empty entries meaning that there is no solution. The entries are as in \autoref{L:RepsSO3TrilinearPreLim}.
\end{Lemma}

\begin{proof}
This is \cite[Theorem 12]{ThCh-trilinear}. See also \cite[Pages 2-4]{Ng-trilinear-moduli} for explicit matrix forms.	
\end{proof}

\begin{Example}\label{E:RepsSO3Trilinear}
To exemplify how to read the table \autoref{Eq:RepsSO3Table}, let us consider 
the column $(3,3,3)$. The only nonzero possibility is that all 
three ternary cubics are of type $xyz=0$. We thus get
\begin{gather*}
\ten_{x}(\bx)=
\begin{psmallmatrix}
t_{111}x_{1}+t_{211}x_{2}+t_{311}x_{3} & 
t_{112}x_{1}+t_{212}x_{2}+t_{312}x_{3} & 
t_{113}x_{1}+t_{213}x_{2}+t_{313}x_{3}
\\
t_{121}x_{1}+t_{221}x_{2}+t_{321}x_{3} & 
t_{122}x_{1}+t_{222}x_{2}+t_{322}x_{3} & 
t_{123}x_{1}+t_{223}x_{2}+t_{323}x_{3}
\\
t_{131}x_{1}+t_{231}x_{2}+t_{331}x_{3} & 
t_{132}x_{1}+t_{232}x_{2}+t_{332}x_{3} & 
t_{133}x_{1}+t_{233}x_{2}+t_{333}x_{3}
\end{psmallmatrix}
\congruence
\begin{psmallmatrix}
x_{1} & 0 & 0
\\
0 & x_{2} & 0
\\
0 & 0 & x_{3}
\end{psmallmatrix}
.
\end{gather*}
That is, we can assume that $t_{111}=t_{222}=t_{333}=1$ and $t_{hij}=0$ otherwise.
\end{Example}

\begin{Lemma}\label{L:RepsSO3TrilinearTwo}
From the cases listed in \autoref{Eq:RepsSO3Table} 
precisely ($(0,0,0)$ and projective type 10) can be used to define an $\son$ tuple up to $\congruence$.
\end{Lemma}

\begin{proof}
Firstly, up to $\congruence$, we have
\begin{gather*}
\mat=
\begin{psmallmatrix}
1 & 0 & 0
\\
0 & 0 & 1
\\
0 & \qpar^{2} & 0
\end{psmallmatrix},
\end{gather*}
by \autoref{E:RepSO3List}. The trilinear form that gives a solution is
\begin{gather*}
t_{123}=-1,\quad
t_{132}=1,\quad
t_{213}=-1,\quad
t_{231}=1,\quad
t_{321}=-1,\quad
t_{312}=1,
\end{gather*}
where we only show the nonzero entries.
That this trilinear form works is a direct calculation.
This is ($(0,0,0)$ and projective type 10) or the \emph{Veronese cuboid}.
All other nonzero cases in \autoref{Eq:RepsSO3Table}
can be directly ruled out.
Since the trilinear form cannot be zero due to the H=I relation, the proof completes.
\end{proof}

\begin{Remark}\label{R:RepSO3EndOfTheLine}
The analog of \autoref{L:RepsSO3Trilinear} for higher forms that would be relevant for \autoref{T:RepsSO3Main}, i.e. 
$(p,p,p)$ trilinear forms with $p\geq 4$, seems to 
be not trackable. In fact, this problem for general $p$ is very difficult, see e.g. \cite{BeSe-matrix-problems}.

\cite{Th-trilinear} has some results regarding 
$(p,p,2p-2)$ trilinear forms, but these are not relevant for $\son$.
For certain subclasses of trilinear forms a bit more can be said, see for example \cite{CoHe-alternating-trilinear}.

The paper \cite{Ng-trilinear-moduli} studies trilinear forms from a 
geometric invariant theory point of view. \cite[Proposition 5]{Ng-trilinear-moduli} gives 
a numerical condition for the stability under GIT quotients of
$(p,q,r)$ trilinear forms. Another geometric treatment is given in \cite{Ng-trilinear-2}, but for $(3,3,4)$ trilinear forms; in particular, the moduli space of such forms 
is related to the moduli space of unordered set of six points in the plane, or dually, six lines. The double
cover of the plane branched along the six lines is a K3 surface, and interesting 
geometry appears. The analog for the $(p,p,p)$ trilinear forms with $p\geq 4$ relevant for this paper appears to be out of reach.
\end{Remark}


\section{On the complexity of the classification problems}\label{S:Complexity}


In this section $\qpar$ is allowed to be any nonzero complex number. It will play the role of a parameter. 

The \emph{rank one classification problem} for 
web categories, say $\web[{\sln}]$, $\web[{\gln}]$ or $\web[{\son}]$, is the classification of 
rank one simple transitive $2$-representations of such categories for all $\qpar\in\C$ at once.
Here classification should be read in the sense of \autoref{C:RepsClassification}.

\begin{Remark}\label{R:ComplexityClassification}
This is again not meant as a definition. In particular, the below 
are \emph{sketchy statements} with \emph{sketchy proofs}. We however hope that we are convincing enough 
so that the reader believes that making these precise (in the sense of complexity theory) 
is not difficult. \changed{We think that making this section precise by properly addressing the complexity questions outlined below is an interesting problem, e.g. is there some finite-tame-2-wild trichotomy for 2-representations?}
\end{Remark}

\changed{$\web[{\sln}]$, and also $\web[{\gln}]$ (or $\qeb$),} is very close to 
be the free pivotal category generated by a bilinear form:

\begin{Proposition}\label{P:ComplexitySL2}
The rank one classification 
problem for $\web[{\sln}]$ implies the 
classification of bilinear forms. Similarly, 
The rank one classification 
problem for $\web[{\gln}]$ \changed{(or $\qeb$)} implies the 
classification of bilinear forms as well.
\end{Proposition}

\begin{proof}
We start by pointing out that all the statements in 
\autoref{SS:RepsSL2Proof} until, and including, 
\autoref{L:RepsSL2MatrixAll} work even if $\qpar$ 
is a nontrivial root of unity. Moreover, \autoref{L:RepsSL2MatrixEquiClasses} 
also holds, but needs to be adjusted as in \cite[Theorem 2.3]{EtOs-module-categories-slq2}.

We assume now that the rank one classification 
problem for $\web[{\sln}]$ is solved. By the above mentioned lemmas
we can associate $\mat\in\matgln$ to a $2$-representation
$\mrep_{\mat}^{n}$ for some $\web[{\sln}]$ by choosing $\qpar$
appropriately. To see this, we point out that the 
relation
\begin{gather*}
\begin{tikzpicture}[anchorbase,scale=1]
\draw[usual] (0,0)to[out=90,in=180] (0.5,0.5)to[out=0,in=90] (1,0);
\draw[usual] (0,0)to[out=270,in=180] (0.5,-0.5)to[out=0,in=270] (1,0);
\end{tikzpicture}
=-\qnum[2]
\end{gather*}
does not give any restriction on the appearing bilinear form 
if we are allowed to vary $\qpar$.
This can be done since $\trace(\mat^{T}\mat^{-1})\in\C$ is some value and we can solve 
$\trace(\mat^{T}\mat^{-1})=-\qnum[2]$ for $\qpar\in\C\setminus\{0\}$.
Thus, we obtain the classification of $\mat\in\matgln$ up to orthogonal congruence (recall that orthogonal congruence is the congruence that preserves the trace). This problem for Hermitian matrices, by 
\cite[Corollary 2.3]{Ho-orthogonal-congruence} and \cite[Theorem 11]{Ri-bilinear}, is equivalent to the classification of nondegenerate bilinear form. The latter 
is then equivalent ot the classification of all bilinear forms, 
as shown in \cite[Unique theorem in Section 1]{Ga-degenerate-bilinear}.

The case of $\web[{\gln}]$ \changed{(or $\qeb$)} can be proven similarly and is omitted.
\end{proof}

We do not know how to deal with the $H=I$ relation, so let us ignore it. Precisely, let $\web[{\son}]^{\prime}$ be the same as $\web[{\son}]$ but without imposing the $H=I$ relation.
The category $\web[{\son}]^{\prime}$ is close to be the free 
pivotal category generated by a trilinear form:

\begin{Proposition}\label{P:ComplexitySO3}
The rank one classification problem 
for $\web[{\son}]^{\prime}$ implies the 
classification of trilinear forms.
\end{Proposition}

\begin{proof}
The proof strategy and arguments are almost the same as in 
the proof of \autoref{P:ComplexitySL2}, so let us only focus on the differences.

We want to argue that the relations
\begin{gather*}
\begin{tikzpicture}[anchorbase,scale=1]
\draw[usual] (0,0)to[out=90,in=180] (0.5,0.5)to[out=0,in=90] (1,0);
\draw[usual] (0,0)to[out=270,in=180] (0.5,-0.5)to[out=0,in=270] (1,0);
\end{tikzpicture}
=\qnum[3]
=\qpar^{2}+1+\qpar^{-2}
,\quad
\begin{tikzpicture}[anchorbase,scale=1]
\draw[usual] (0,0)to[out=90,in=180] (0.5,0.5);
\draw[usual] (1,0)to[out=90,in=0] (0.5,0.5);
\draw[usual] (0.5,0)to (0.5,0.5);
\draw[usual] (1,-0.5) to (1,0);
\draw[usual] (0,0)to[out=270,in=180] (0.25,-0.5)to[out=0,in=270] (0.5,0);
\end{tikzpicture}
=0
,
\end{gather*}
will not restrict the choice of trilinear form.

Similarly as in the proof of \autoref{P:ComplexitySL2} we can vary
$\qpar$, eliminate the circle relation and we can assume that 
$\mat\in\matgln$ (encoding the bilinear form) is arbitrary. The bitri evaluation thus does not restrict the appearing trilinear form because we can just chose the bilinear form accordingly. \changed{Although this is difficult in practice, this follows from a simple parameter count. Namely, the trilinear form has $n^{3}$ parameters, and so its kernel has $n^{3}-1$ parameters, while the bilinear form has $n^{2}$. (Note that this count does not work for $n=1$, but classifying $(1,1,1)$-trilinear forms is trivial.)}

The remaining steps work as at the 
end of the proof of \autoref{P:ComplexitySL2} (the reduction from nondegenerate trilinear to general trilinear forms follows by copying the proof of \cite[Unique theorem in Section 1]{Ga-degenerate-bilinear}).
\end{proof}

Note that all of our rank one classification 
problems have an associated $\C$-vector space,
i.e. the image of the 
generating object. Let $n\in\N$ denote the dimension of this space.

In analogy with matrix classification problems, 
we call a rank one classification 
problem \emph{finite} if there are only finite many equivalence 
classes of rank one simple transitive $2$-representations for every fixed $n\in\N$. Similarly, such a problem is \emph{tame} if there 
is at most a one-parameter family of equivalences classes per $n$.
We call such a problem \emph{2-wild} 
(alternatively, \emph{wilder than wild}) if it is strictly 
more difficult than any wild problem in the sense that solving it 
solves all wild problems, but not vice versa. (Recall that a classification problem 
is called \emph{wild} if it contains the classification of 
indecomposables for any finite dimensional algebra.)

\begin{Theorem}\label{T:ComplexityDifficult}
The rank one classification problem for $\web[{\son}]^{\prime}$ is
2-wild.
\end{Theorem}

\begin{proof}
This follows from \autoref{P:ComplexitySO3} and \cite[Theorem 1.1]{BeSe-matrix-problems}.
\end{proof}

The above, together with the easy to obtain solution of the rank one 
classification problem for $\web[{\sln[1]}]$, implies the following.
The rank one classification problem...
\begin{enumerate}[label=(\roman*)]

\item ...for $\web[{\sln[1]}]$ is finite.

\item ...for $\web[{\sln}]$/$\web[{\gln}]$ is tame.

\item ...for $\web[{\son}]^{\prime}$ is 2-wild.

\end{enumerate}

\begin{Remark}\label{R:ComplexityDifficult}
\leavevmode	
\begin{enumerate}

\item In the representation theory of finite dimensional algebras 
there is the notion of \emph{finite, tame and wild representation type}. The above is inspired from these 
notions.

\item Note that the categorical version of wild, that we called 
2-wild, is strictly more difficult than any wild problem. 
In this sense one can say that categorical representation theory is more difficult than classical representation theory. However, 
the main caveat is that we are discussing $\web[{\son}]^{\prime}$
and not $\web[{\son}]$ itself.

\end{enumerate}
We think it is an interesting question whether the 
rank one classification problem for $\web[{\son}]$ 
(and with it probably for most other web categories)
is wilder than wild.

Optimally, we would like to write the rank one classification problem for $\web[{\sln[1]}]$ is finite, 
for $\web[{\sln[2]}]$ it is tame and for $\web[{\sln[3]}]$ 
it is 2-wild. (In order, no form appear, bilinear forms appear, and trilinear forms appear.) We however were not able to verify this because of the so-called square relation for $\web[{\sln[3]}]$.
\end{Remark}


\newcommand{\etalchar}[1]{$^{#1}$}

\end{document}